\documentclass[preprint,10pt]{elsarticle}
\usepackage{helvet}
\usepackage[utf8]{inputenc}
\usepackage{hyperref} 

\usepackage{amsmath,amsfonts,amssymb}
\usepackage{mathrsfs}
\usepackage{lineno}
\usepackage{bm}
\usepackage{stmaryrd}
\usepackage{graphicx}
\usepackage{nicefrac}
\usepackage{multirow}

\newcommand{\cT}{ \mathcal{T} }

\usepackage{algorithm2e}
\RestyleAlgo{ruled}
\graphicspath{ {./Figures/} }

\usepackage{caption}
\usepackage{subcaption}
\usepackage{xcolor}
\usepackage{epstopdf}

\usepackage{comment}

\usepackage{scalerel}
\DeclareMathOperator*{\A}{\scalerel*{A}{\sum}}

\usepackage{tikz}
\usetikzlibrary{positioning,backgrounds,fit}
\usetikzlibrary{calc}
\usepackage{scalerel}

\DeclareMathAlphabet{\mathantt}{OT1}{antt}{li}{it}
\DeclareMathAlphabet{\mathpzc}{OT1}{pzc}{m}{it}

\begin{document}

\begin{frontmatter}

\title{Three-dimensional modelling of drag anchor penetration using the material point method}

\author[dur]{Robert E. Bird}
\author[dur]{William M. Coombs\corref{cor1}}\ead{w.m.coombs@durham.ac.uk}
\author[dundee]{Michael J. Brown}
\author[dur]{Charles E. Augarde}
\author[dundee]{Yaseen U. Sharif}
\author[dur]{Giuliano Pretti}
\author[bgs]{Catriona Macdonald}
\author[bgs]{Duncan Stevens}
\author[arup,bgs]{Gareth Carter}

\cortext[cor1]{Corresponding author}

\address[dur]{Department of Engineering, Durham University, South Road, Durham, DH1 3LE, UK}
\address[dundee]{School of Science and Engineering, University of Dundee, Fulton Building, Dundee, DD1 4HN, Scotland, UK}
\address[arup]{Arup, 10 George St, Edinburgh, EH2 2PF, Scotland, UK}
\address[bgs]{British Geological Survey, Currie, Edinburgh, EH14 4BA, Scotland, UK}

\begin{abstract}
\noindent Drag embedment anchors are a key threat to buried subsea linear infrastructure, such as power/data cables and pipelines. For cables, selecting a burial depth is a compromise between protecting the cable from anchor strike and the increased cost of deeper installation. This presents an efficient large deformation, elasto-plastic Material Point Method-based soil-structure interaction predictive tool for the estimation of anchor penetration based on Cone Penetration Test (CPT) site investigation data. The tool builds on earlier work by the authors supplemented by three developments: modelling assemblies of rigid bodies (necessary for articulated anchors), a partitioned domain approach to enable accurate and efficient modelling of long anchor pulls and an improved means of modelling rotational inertia. The tool is validated against scaled physical tests conducted in a geotechnical centrifuge on sands with a range of relative densities with good agreement across the tested conditions. Numerical simulations identify key issues with the UK Cable Burial Risk Assessment (CBRA) approach for estimating anchor penetration and reveal the potentially non-conservatism of the CBRA framework for sandy seabeds. The numerical model enables site-specific anchor-penetration assessment along cable routes and can be used to evaluate the performance of different anchor designs and sizes in varied soil conditions.

\end{abstract}


\end{frontmatter}

\section{Introduction}
Subsea cables provide vital arteries for data and power transfer worldwide. In shallow waters these cables are buried below the seabed surface in order to protect them from accidental interaction with anchors/fishing gear that, based on 2024 data, cause around 70\% of the reported cable failures \cite{cables2024}. For power cables the burial depth (or \textit{depth of lowering}) along a cable route can be specified using the UK Carbon Trust's Cable Burial Risk Assessment (CBRA, \cite{CBRA}) framework. This approach relies on anchor-seabed penetration factors used to estimate the penetration potential of an anchor as a multiplier of its fluke length. However, the CBRA framework provides factors for just two broad soil groups: a factor of $1$ for sands/stiff clays and a factor of $3$-$5$ for soft clays, and  acknowledges that \emph{"industry as a whole would benefit from further research into anchor penetration in a range of seabed types"} \cite{CBRA}. The primary motivation for the research presented in this paper is to provide more nuanced seabed penetration factors obtained from simulations of anchor penetration using a numerical model calibrated using routine site investigation data. To achieve this here we extend the large deformation Material Point Method (MPM) soil-structure interaction framework of Bird \emph {et al.} \cite{bird2024dynamic} to include new features. 

A variety of methods exist to determine the kinematics of drag anchors in soils, which can broadly be categorised as analytical, experimental, and numerical. One of the first articles to present  analytical analysis for the drag anchor trajectory was Neubecker and Randolph \cite{neubecker1996static,neubecker1996kinematic}, where an energy minimisation process was performed to determine the incremental direction of the anchor (incremental translation or rotation). Other notable articles include Stewart \emph{et al.} \cite{stewart1992drag} where an analytical force and moment equilibrium is used to determine the increment motion of the anchor. Additionally the hybrid numerical-analytical method of O'Neil \emph{et al.} \cite{o2003drag} which uses a finite element analysis to first define the load yield locus of the admissible anchor loads in a elasto-plastic domain, through incremental motions of the fluke (horizontal, vertical and rotational). The incremental motion of the anchor is then obtained through force, and moment, equilibrium, of the anchor with the chain tension, which locates the current load point on the yield locus and thus the corresponding incremental motion. Other, more recent analytical methods, include Liu \emph{et al.} \cite{LIU201233} who derive a closed form kinematic expression for the anchor where the depth at time $t$ is defined by the reverse catenary of the chain and Ren \emph{et al.} \cite{REN2022110699} where a force-moment equilibrium, using similar methods to \cite{neubecker1996static}, is performed with a forward-Euler time integration to determine the anchor motion. Analytical models have the benefit that they are very fast compared to numerical methods however, as discussed in \cite{pretti2022review,pretti2024continuum}, whilst the simplicity of analytical methods makes them attractive, if parameters for the particular anchor and soil conditions are unknown then these formulations cannot accurately, or reliably, make predictions for this non-linear (materially and geometrically) dynamic problem involving multiple bodies.

Numerical methods allow for flexibility in anchor geometry, material behaviour and anchor deployment conditions. Possibly the largest numerical challenge in modelling drag anchors is the very large deformation of the soil, the seemingly most common approach is to use the explicit-in-time Coupled-Eulerian-Lagrangian (CEL) formulation in Abaqus FEA \cite{abaqus2011}, for example see the non-exhaustive list \cite{liu2014numerical,grabe2015numerical,dao2024numerical,DOU2018199,zhao2016numerical,zhao2016efficient,ZHANG2023105518}. The studies cover a wide range of topics, such as investigating the effect of installation velocities on the general kinematics of the anchor for varying soil conditions, Liu \emph{et al.} \cite{liu2014numerical}, or, how an anchor's geometry effects its trajectory \cite{DOU2018199} and also, more recently, real offshore site prediction for the performance of a range of drag anchors in various sands, Dou \emph{et al.} \cite{dao2024numerical}. However, an issue with numerical modelling of anchors is: (i) the cost of a simulation is expensive (anchor modelling is inherently 3D) and (ii) it can produce results that appear accurate but are not fully representative of reality, confidence in the results can only be achieved through validation. Currently the works in which numerical validation is performed is limited, and often the validation itself is limited to only a single anchor type and/or a single material. However, there is experimental data available to validate against for drag anchors, O'Neil and Randolph \cite{o2001modelling} who model anchor penetration in a drum centrifuge, Moore and co-workers \cite{Moore2017,MOORE2021109411} also use a drum centrifuge for determining anchor penetration characteristics, Liu \emph{et al.} \cite{LIU2010434} who experiment on drag anchor designs with various geometries in a $1$g setup, Yan \emph{et al.} \cite{YAN2015529} who experiment on rock berms, and Sharif et al \cite{Yaseen2024} who perform centrifuge, and $1$g, tests for a single anchor but for a large range in sand properties. Informed by the recent reviews of \cite{pretti2022review,pretti2024continuum}, this work validates the proposed modelling approach against experimental drag anchor pulls for a wide range in sand relative densities. The experimental data is provided by Sharif et al \cite{Yaseen2024} for rate independent fully drained conditions, see the work of Brown et al. \cite{Brown2026rate} where rate effects are quantified. additionally this work introduces a new method that decouples simulation time per step from the total drag length. 

The numerical model described and demonstrated in our previous work \cite{Bird2024,bird2024dynamic} comprises a single rigid body interacting with a soft highly deformable material (e.g. the seabed). The latter is represented using the Material Point Method (MPM), which can cope with large deformations in solids much better than the standard Finite Element Method (see \cite{Solowski2021} and \cite{Vaucorbeil2020} for review articles). Key advances described in \cite{Bird2024,bird2024dynamic} necessary for the model to function include a new and efficient contact algorithm which provides a general and consistent representation of the extent of the deformable body without introducing boundary representation in the MPM.

The ability to model multiple rigid bodies within a single simulation becomes essential when aiming to accurately represent offshore geotechnical assemblies that incorporate hinges or joints. Examples include the skid plates of seabed cable and pipeline ploughs, the interconnected tracks of subsea Remotely Operated Vehicles (ROVs), and the fluke–shank configuration of drag embedment anchors. To address this need, here we extend our previous model to support assemblies of multiple rigid bodies, enabling the simulation of realistic articulating anchors. In doing this, enhancements have also been made to the representation of rotational inertia for these bodies. A third key development is the introduction of a computationally efficient numerical scheme for simulating realistic anchor pulls. For such simulations to be practically useful, they must produce data over pull lengths many times greater than the anchor’s own dimensions, as drag anchors often require substantial horizontal displacement to embed fully, depending on both anchor geometry and soil conditions. Our solution is to partition the soft material problem domain as will be explained below.

Having described the fundamental components and the three developments noted above, the predictive tool is demonstrated starting with material calibration using Cone Penetration Tests (CPTs) against experimental data. The resulting parameters are then used as inputs for the drag anchor prediction depths, the results of which are validated against penetration depths obtained from scaled physical tests conducted in a geotechnical centrifuge at the University of Dundee \cite{Yaseen2024}. Of note is that only four mechanical parameters are used in the soil constitutive models (two elastic, two plastic), which provides a good balance between practical use of the model and accuracy of the delivered results. The code is used to investigate the impact of anchor geometry and mass on the penetration depth and capacity in different density sands. The paper closes with a discussion on the implications of these results for the CBRA framework and the performance of drag anchors more generally. 

\section{Model components}\label{sec: Material point method}

While the developments of the MPM to be described below could be used to model a wide range of complex soil-structure interaction problems, here we will focus on the anchor pull problem. The numerical model comprises (1) a model of the seabed soil and (2) a model of the anchor. The former uses the MPM and the approach is based on the three-dimensional implicit dynamic MPM formulation in \cite{bird2024dynamic} details of which are only summarised here. This is a dynamic extension of the quasi-static open-source code AMPLE (A Material Point Learning Environment) code \cite{coombs2020aample} that has its origins in the work of Charlton \emph{et al.} \cite{charlton2017igimp} and Coombs \emph{et al.} \cite{coombs2020on}.  The particular variation of the MPM used here is the General Interpolation Material Point Method (GIMPM) due to its ability to reduce cell crossing instabilities under large deformations. As discussed in \cite{Bird2024,bird2024dynamic} the General Interpolation Material Points (GIMPs) have an associated cuboid domain, the vertices of which are used for contact detection with rigid body surfaces and to apply contact boundary conditions, which are enforced with a penalty method \cite{wriggers2006computational}. This enables an explicit contact formulation based on the GIMP domain and rigid body geometries and resolves the usual lack of a clear definition of the boundary associated with the MPM \cite{acosta2021development,pretti2024continuum}. 

The continuum formulation adopted for the soil in the paper is also the same as in Bird \emph{et al.} \cite{bird2024dynamic}. Full details of the large deformation elasto-plastic continuum mechanics approach can be found in Charlton \emph{et al.} \cite{charlton2017igimp} and Coombs and Augarde \cite{coombs2020aample}. The  approach is based on a Hencky material assumption that adopts a linear relationship between logarithmic (or Hencky) strains and the Kirchhoff stress. This combination allows the format of infinitesimal strain constitutive models to be recovered and used within a large deformation (i.e. finite strain) boundary value solution without modification of the stress update algorithm (see \cite{SouzaNeto} for details). This approach is widely used across finite element and MPMs, see for example \cite{charlton2017igimp,coombs2020aample,coombs2025aggregated,coombs2022ghost,bird2024dynamic,Bird2024,coombs2020on,Pretti2024}. When the GIMPM and rigid body are in dynamic equilibrium, the weak statement of equilibrium is discretised in time using the Newmark method with the parameters $\gamma = 1$ and $\beta = 1/2$. A FLIP node-to-point mapping is performed \cite{brackbill1986flip} at the end of the time step to update the material point velocities, the accelerations do not need updating due to the values of $\gamma$ and $\beta$.

The non mesh-matching nature of the MPM means that it can suffer from the \emph{small cut} issue and associated loss of coercivity of the linear system of equations. In the MPM the problem occurs when nodes of the background mesh have very small mass or stiffness contributions but significant residual force values, which results in very large non-physical accelerations or displacements at the boundary of the physical body.  This issue was recognised by Sulsky~\emph{et al.} \cite{Sulsky1995} and motivated the widely used Modified Update Stress Last (MUSL) explicit algorithm. However, there are currently only two general approaches for resolving this issue for implicit MPMs: ghost stabilisation \cite{Burman2010,coombs2022ghost} and mesh aggregation \cite{coombs2025aggregated}. Consistent with Bird~\emph{et al.} \cite{bird2024dynamic}, in this paper ghost stabilisation is adopted, which adds a penalty term to the mass and/or stiffness matrix that constrains the gradient of the solution across faces of the background mesh near the boundary of the deformable body.

\subsection{Numerical framework for soil-structure interaction}
For completeness the weak forms for the GIMPM and rigid body are presented here, it is a summary of the formulation from Bird \emph{et al.} \cite{Bird2024} where further detail on the derivation, implementation and validation can be found, particularly for the contact. The discrete weak form of the deformable body in an updated Lagrangian frame is,
\begin{align}
   \A_{ \forall p} \Bigl([\nabla_x S_{vp}]^{T}\{\sigma_p\} V_p -[S_{vp}]^{T}\{b\} V_p +[S_{vp}]^{T}\{\dot{v}_p\} m_p \Bigr) &\nonumber\\ \label{eqn:MPMweak} 
   - \A_{p\in P_c}\left(\{\delta f_{N,vp} ^{\varphi}\} + \{\delta f_{T,vp} ^{\varphi}\}\right) &= \{0\},
\end{align}
where $p$ in general denotes quantities associated with a material point, $V_p$ and $m_p$ are respectively the material point volume and mass. $\{\cdot\}$ denotes a vector quantity with $\{\sigma_p\}$ as Cauchy stress and $\{b\}$ the body force. The GIMPM basis functions are expressed as $[S_{vp}]$, where $[\cdot]$ is a matrix, they depend on the convolution of the background mesh basis functions and the characteristic function associated with the material point, \cite{bardenhagen2004generalized}. The set of all material points in contact is defined $P_c$, where the normal and tangential contact forces acting on the material points are respectively, $\{\delta f_{N,vp} ^{\varphi}\}$ and $\{\delta f_{T,vp} ^{\varphi}\}$. 

The engineering machinery used in subsea process often have multiple separate components which interact with the seabed. For the purposes of modelling, each component is represented by a separate rigid body, with multiple bodies connected together with a finite element truss frame. The truss nodes model the kinematics of the rigid bodies with a component's centre of mass set at a truss node. Additionally the machinery is sufficiently large that rotational inertia influences the kinematics. Therefore rotational inertia terms are included to consider the resistance to angular acceleration about a rigid body's centre of mass. The surface of the rigid body is constructed from a continuous set of triangles $\cT$, where $\cT_c\subseteq\cT$ is the set of triangles in contact. As the kinematics are model by truss elements, the contact forces on the surface of the rigid body must be mapped from the contact triangles surface of the rigid body, to the nodes of the truss element $e$, this is achieved with the mapping $\Xi^e[\cdot]$, see Bird \emph{et al.} \cite{bird2024dynamic}. The discrete weak form for the truss frame is therefore
\begin{equation}\label{equ: truss weak form contact}
\begin{split}
 &\A_{e\in E}\int_{\varrho_t(e)}\left([\zeta]^T \{f^I\}+\rho^e[\zeta]^T\{\dot{v}^{e}\}+ [\zeta]^T
 [{\partial\theta}/{\partial x_e}]^T[I(x_e)]
 \{\ddot{\theta}(x_e)\}\right) \text{d}V
  = \\
   &\A_{e\in E}\int_{\varrho_t(e)} [\zeta]^T\{b\} \text{d}V +
   \A_{e\in E}\Xi^e\left[\int_{\varrho_t(\cT_c)}\left(\{\delta f_N^{\varrho}\} + \{\delta f_T^{\varrho}\}\right)\text{d}V\right]  
   \end{split}
\end{equation}
$\varrho_t(\cdot)$ is the motion of the rigid body, a truss element is defined $e$ and belongs to the set of all elements $E$. Additionally, $[\zeta]$ is the matrix of shape functions for the truss element, $\{f^I\}$ is the internal force vector, $\rho^e$ is the material density in $e$ and $\{\dot{v}_{e}\}$ is acceleration. The third term in the left hand side of Equation \eqref{equ: truss weak form contact} corresponds to the rotational inertia of the truss element in three dimensions, where the angle $\theta$, and its derivatives with respect to time, are a function of the truss element's nodal positions in space $x_e$. Additionally the rotational inertia matrix $[I(x_e)]$ will evolve as the truss element rotates. The description of the application, first variations and linearisation of this term about one axis, as only one axis is necessary for modelling the anchor for CBRAs, are are presented in Section \ref{sec: Multiple rigid body interaction}.

The deformable and rigid body systems, respectively Equations \eqref{eqn:MPMweak} and \eqref{equ: truss weak form contact} are coupled together through the contact terms. Each time step is solved implicitly using the Newton-Raphson (NR) method combined with a basic line search which chooses the lower residual corresponding to the line search length $\alpha = \{0.1,1.0\}$ \cite{}. It is was found that a line search was required to prevent contact points repeatedly coming in and out of contact during a load step \cite{dennis1996numerical}.

The linear system for each NR iteration takes the form 
\begin{equation}\label{equ: NR incremental solution}
\underbrace{
  \left\{
    \begin{array}{c}
      \{\Delta u_{v}\} \\
      \{\Delta u_{b}\}
    \end{array}
  \right\}
}_{\{\Delta u\}}
    =
   \underbrace{ \left[
    \begin{array}{cc}
  {[A]} & {[B]} \\
  {[C]} & {[D]}
\end{array}
    \right]^{-1}}
    _{[K]^{-1}}
    \underbrace{ \left\{
      \begin{array}{c}
         \{f_{v}\} \\
           \{f_{b}\}
    \end{array}\right\}}_{\{f\}}
\end{equation}
$\Delta u$ is the incremental solution update and $f$ is the out-of-balance forces of the system. $[A]$ and $[B]$ are matrices from the linearisation of Equation \eqref{eqn:MPMweak}, whilst $[C]$ and $[D]$ are from the linearisation of Equation \eqref{equ: truss weak form contact}. The subscript $v$ is used for terms in Equation \eqref{equ: NR incremental solution} to describe the nodes (vertices) of the background grid on which the GIMPM is solved at each time step, the subscript $b$ is used to describe the DoF of the truss element frame associated with the rigid bodies. For all considered problems the boundary conditions were background mesh confirming, and therefore the same techniques as finite elements could be used in the GIMPM framework \cite{coombs2020aample} with Dirichlet boundary conditions applied directly to the vertices $v$ of the mesh. Any rigid body constraints are also imposed directly. 

\subsection{Modelling multiple rigid bodies}\label{sec: Multiple rigid body interaction}
Having given details of the general components of the numerical model, here we describe in more detail the new features added to enable the anchor pull problem to be efficiently modelled. The contact formulation used here matches that in the previous model \cite{bird2024dynamic} where the surface of a rigid body is formed from a set of flat triangles $\Upsilon\in\cT$, which provide the surface component of the \textit{point-to-surface} contact. The corners of all the GIMP domains form the deformable domain's contact surface $p\in{P_c}$, and provide the point component of the \textit{point-to-surface} contact. A gap function, $g_N$, obtained from a Closest Point Projection (CPP) algorithm is used to determine whether contact occurs and a penalty law minimises the normal interpenetration. The CPP projects the corner of the GIMP $\bm{x}$ onto a point $\bm{x}^\prime$ on one of the triangles representing the rigid body. The set of all triangles in contact is defined $\Upsilon\in\cT_c$.

By definition the rigid bodies cannot deform, they can only rotate and translate as a function of the rigid body's kinematics. Therefore a definition of the rigid bodies' kinematics is required which also defines the positions of the surfaces, necessary for contact. Each triangle has its own local coordinate system $\xi^\alpha$, where $\alpha\in[1,2]$, and the position of points on the triangle can be expressed using the usual linear triangle basis functions, $N_i(\xi^\alpha)$, that is
\begin{equation}\label{eq: tri position}
    \bm{x}^\prime = \sum_{i = 1}^{3}\sum_{\alpha=1}^2N_i(\xi^\alpha)\bm{x}_{\Upsilon,i}^\prime,
\end{equation}
where bold notation is used to define a vector, $i$ refers to the basis number and $\bm{x}_{\Upsilon,i}^\prime\in\mathbb{R}^3$ are the positions of the triangle's vertices in three dimensional space. 

\begin{figure}[ht!]
    \centering
    \begin{subfigure}[t]{0.35\textwidth}
        \centering
       \includegraphics[width=0.9\textwidth]{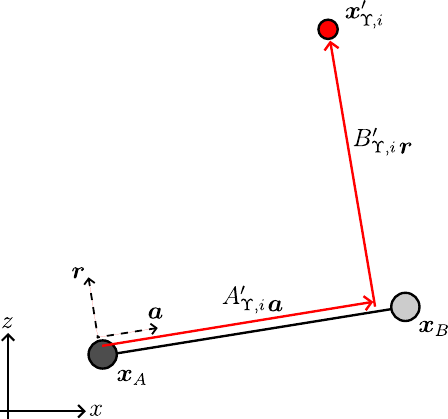}
        \caption{} \label{fig:truss frame}
    \end{subfigure}%
    ~ 
    \begin{subfigure}[t]{0.62\textwidth}
        \centering
      \includegraphics[width=0.9\textwidth]{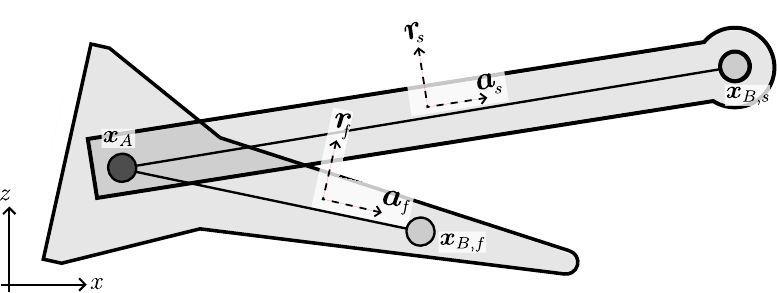}
        \caption{} \label{fig:truss frame with anchor}
    \end{subfigure}
    \caption{Rigid body: (a) is a single line element with its normal $\bm{a}$ and tangent $\bm{r}$ being used to define the point $\bm{x}_{\Upsilon,i}^\prime$. (b) shows how two line elements are used to define the coordinates of the fluke and shank of an anchor.}
\end{figure}

A truss frame finite element method is used to join together the multiple rigid bodies, all of which are in contact with the deformable body. This method was validated for a single rigid body in \cite{bird2024dynamic} where the position of its surface is a function of the nodal positions of a single truss element and is independent of the element's deformation\footnote{The stiffness of the element is sufficiently high that only arbitrarily small deformations occur.}. An example is shown in Figure \ref{fig:truss frame} with a truss element which has nodes $A$ and $B$ at positions $\bm{x}_A$ and $\bm{x}_B$ respectively. From the node positions of the truss the axial and radial directions are also respectively defined,  
\begin{equation}
    \bm{a} = (\bm{x}_B-\bm{x}_A)/|\bm{x}_B-\bm{x}_A|\quad\text{and} \quad \bm{r} = \bm{R}\cdot\bm{a}.
\end{equation}
where $\bm{R}$ is a rotation matrix of $90^\circ$ about the $y$-axis. To determine, $\bm{x}_{\Upsilon,i}^\prime$, the position of vertex $i$ of triangle $\Upsilon$, necessary for Equation \eqref{eq: tri position}, the following expression is used
\begin{equation}\label{eq: vertex equation}
\bm{x}_{\Upsilon,i}^\prime=\bm{x}_{\Upsilon,i}^\prime(\bm{x}_A,\bm{x}_B) = A^\prime_{\Upsilon,i}\bm{a} + B^\prime_{\Upsilon,i}\bm{r} + \bm{x}_A,
\end{equation}
where $A^\prime_{\Upsilon,i}$ and $B^\prime_{\Upsilon,i}$ are constants that respectively define how far in the axial and radial direction the point $\bm{x}_{\Upsilon,i}^\prime$ is from $\bm{x}_A$. The positions $\bm{x}_A$ and $\bm{x}_B$ are the primary variables of the rigid body which are solved for, they fully define the kinematics of the rigid body with the vertex $\bm{x}_{\Upsilon,i}^\prime$ a function of their position.

For the multiple body case, each rigid body is defined with a truss element, which are connected together permitting rotation of the individual bodies with respect the each other. No interpenetration or contact is permitted hence the surfaces from different rigid bodies are non-intersecting, meaning their interiors are disjoint. As an example, an anchor comprising two components, the fluke $f$ and the shank $s$, is shown in Figure \ref{fig:truss frame with anchor}. Each are modelled as rigid bodies with the surfaces constructed from triangles, $\Upsilon\in\cT_s$ and $\Upsilon\in\cT_f$, where the set of all triangles is defined as $\cT = \cT_s\cup\cT_f$. The trusses in Figure \ref{fig:truss frame with anchor} share the node $A$ and have their own axial and radial directions which are used to define the position of the their respective surfaces. The primary unknowns for this system are the positions $\bm{x}_A$, $\bm{x}_{B,s}$ and $\bm{x}_{B,f}$.

Appropriate representation of a structure's rotational inertia is key in predicting its kinematic behaviour if there is any unconstrained rotation. This is particularly the case for anchors which can undergo large rotations, at various rates, and often can be several metres long and tens of tonnes\footnote{For an example on the range of commercial anchor size see the Sotra product page \cite{anchor_website}.}. Within the framework of the truss frame's translation Degrees of Freedom (DoF) there are two options for modelling the inertia: Option 1 is to place a series of points around the Centre of Mass (CoM) at a correct radius with a distributed mass to represent the rotational inertia, as in \cite{bird2024dynamic}; Option 2 is to define a rotation, and acceleration, as a function of the translational DoF. Option 2 is used here as it is exact but requires more work to determine its linearisation whereas, Option 1 is inexact and could have stability issues if a large number of stiff truss members are used.

\begin{figure}[ht!]
   \centering
       \includegraphics[width = 0.4\textwidth]{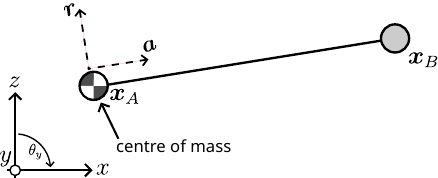}
        \caption{Rotational inertia: An example of a truss element with angular acceleration defined at the centre of mass $\bm{x}_A$.} \label{fig:truss frame with rot}
\end{figure}

Although the approach adopted in this paper is general and can be applied to any rigid body undergoing soil-structure interaction, the subsequent explanation focuses on the specific case of an articulated drag anchor that can rotate about a a single axis. Considering Figure \ref{fig:truss frame with rot}, to include rotational inertia a description of the angular acceleration around the centre of mass is required, this can be achieved using the end points of the line element, $\bm{x}_A$ and $\bm{x}_B$. The points $\bm{x}_A$ and $\bm{x}_B$ define the axial direction of the element. The rotational inertia term which contributes to the rigid body's weak form, and respectively its linearisation, takes the form 
\[
\delta f^\theta = I_{yy}~\delta\theta\ddot{\theta}\quad\text{and}\quad \Delta\delta f^\theta =I_{yy}~\Delta\delta\theta\ddot{\theta} + I_{yy}~\delta\theta\Delta\ddot{\theta}
\]
Using the notation $\left(\cdot\right)^{AB}  =  (\cdot)_B- (\cdot)_A$ the rotation, angular velocity and acceleration can be defined
\begin{equation}
    \theta = \tan^{-1}\left(\frac{{x}^{AB}_3}{{x}^{AB}_1}\right),\quad \dot{\theta} = \frac{\partial\theta}{\partial t}= \frac{\partial\theta}{\partial \bm{x}^{AB}}\cdot\frac{\partial \bm{x}^{AB}}{\partial t}
\end{equation}
and
\begin{align}
    \ddot{\theta} &= \frac{\partial}{\partial t}\left(\frac{\partial\theta}{\partial \bm{x}^{AB}}\cdot\frac{\partial \bm{x}^{AB}}{\partial t}\right)
     = \frac{\partial \bm{x}^{AB}}{\partial t}\cdot\frac{\partial^2\theta}{\partial {\bm{x}^{AB}}^2}\cdot\frac{\partial \bm{x}^{AB}}{\partial t} +
   \frac{\partial\theta}{\partial \bm{x}^{AB}}\cdot\frac{\partial^2 \bm{x}^{AB}}{\partial t^2}\\
   &= \bm{v}^{AB}\cdot\frac{\partial^2\theta}{\partial {\bm{x}^{AB}}^2}\cdot\bm{v}^{AB} +  \frac{\partial\theta}{\partial \bm{x}^{AB}}\cdot\dot{\bm{v}}^{AB}
\end{align}
with the variation, $\delta\theta$, and linearisations $\Delta\ddot{\theta}$ and $\Delta\delta\theta$ provided in \ref{App:Linearisation of theta}.

\begin{figure}[ht!]
    \centering
    \begin{subfigure}[t]{0.4\textwidth}
        \centering
       \includegraphics[width=1\textwidth]{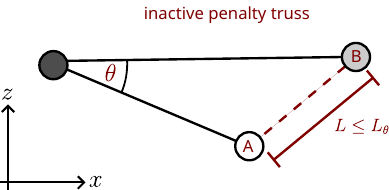}
        \caption{} \label{fig:truss_frame_new_rotation_1}
    \end{subfigure}%
    ~ 
    \begin{subfigure}[t]{0.4\textwidth}
        \centering
      \includegraphics[width=1\textwidth]{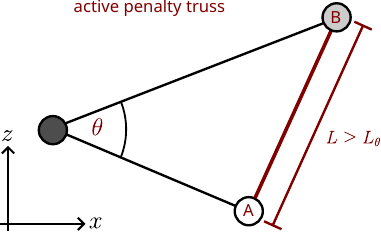}
        \caption{} \label{fig:truss_frame_new_rotation_2}
    \end{subfigure}
    \caption{Limiting opening angle: When the opening angle is less than the limit the penalty spring is inactive (a) but when greater than a limit it is active (b).}
\end{figure}

In practice, it is often also necessary to constrain the relative rotation between the rigid bodies in an assembly. In the case of a drag anchor with a pivoting shank and fluke, the opening angle of the anchor is constrained by the physical geometry of the connection. In this paper a penalty approach is used to introduce these types of constraint.  An example is shown in Figures \ref{fig:truss_frame_new_rotation_1} and \ref{fig:truss_frame_new_rotation_2}, where the angle $\theta$ is limited by including a penalty spring between nodes A and B. To determine when the penalty stiffness should be active the distance between nodes is calculated for the maximum allowable value of $\theta$, this distance is defined $L_\theta$. $L_\theta$ is also defined as the undeformed length of the penalty spring. During a simulation if the distance between nodes A and B is less than $L_\theta$ there is no additional penalty force between nodes A and B (the case shown in Figure~\ref{fig:truss_frame_new_rotation_1}). When the distance is greater than $L_\theta$ and additional penalty force is introduced between nodes A and B to limit the angle $\theta$ (the case shown in Figure~\ref{fig:truss_frame_new_rotation_2}).

\subsection{A partitioned domain approach}\label{sec: Partitioned domain approach}
\begin{figure}[ht!]
    \centering
    \begin{subfigure}[t]{0.49\textwidth}
        \centering
       \includegraphics[width=1\textwidth]{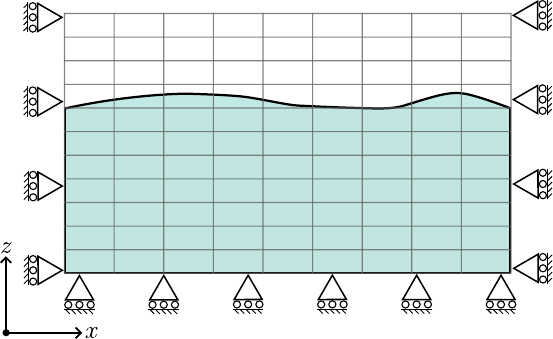}
        \caption{} \label{fig:partdomain stage 1}
    \end{subfigure}%
    ~ 
    \begin{subfigure}[t]{0.49\textwidth}
        \centering
      \includegraphics[width=1\textwidth]{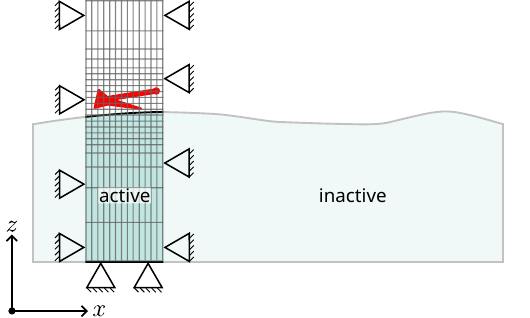}
        \caption{} \label{fig:partdomain stage 2}
    \end{subfigure}
    \caption{Partitioned domain: The mesh and initial material configuration for stage 1 is shown in (a) with the subsequent deformed material, new mesh and boundary conditions, and anchor (in red) for stage 2 shown in (b).}
\end{figure}

The previous sections have described the first two developments in the new model and in this section the third is detailed. With respect to drag anchor pulls, it is never known \emph{a priori} how far an anchor will need to be dragged in order to achieve its ultimate penetration depth. This can lead to two potential problems in numerical modelling: (i) the length of the domain may be insufficient for the anchor to achieve its ultimate depth, resulting in an analysis rerun with a longer domain; or (ii) the length of the domain may be longer than required to obtain the ultimate depth, resulting in wasted computational effort. To mitigate these issues a new \emph{partitioned domain} approach is introduced here. The aim is to make the solve time per load step, for a given refinement, invariant to the total drag length with minimal effect on the anchor trajectory. The method is designed so it is possible to include variation in material properties and surface topology.
In the simulation of an anchor pull there are two distinct stages of loading: firstly the anchor is placed on the seabed and embeds slightly under self-weight. The seabed itself will be under stress due to gravity. Secondly the anchor is translated and the forces on the anchor arise from contact between the rigid bodies and the soft deformable domain. Here we model this by undertaking a two stage analysis. We define a large deformable material domain using the MPM but use it only in one of the two analysis stages.
\begin{description}
    \item[Stage 1:] A quasi-static simulation to determine the stress state and material deformation under gravitational load, Figure \ref{fig:partdomain stage 1}; and
    \item[Stage 2:] A dynamic simulation of the anchor being pulled through the deformed material Figure \ref{fig:partdomain stage 2}.
\end{description}
Stage 1 is solved for the entire domain, which ensures that the initial stress state, and displacement, are correctly determined for features such as, variations in material properties and surface topology. Stage 1 calculates the initial material state that is used as an input to Stage 2. In Stage 2 a domain around the anchor is defined which is a partition of the original larger domain. Material points which are in the partitioned domain are set to be active whilst material points outside the domain are set to be inactive and are fixed. As the anchor moves through the domain the partitioned domain moves with it. During movement, material points will leave the partition becoming inactive, whereas other material points, with their initial state calculated in Stage 1, will enter the partitioned domain and become active in the analysis. The overall simulation cost is a function of the mesh used and the GIMP resolution inside the partitioned domain, and hence a greater resolution around the anchor can be achieved for the same computational cost compared to considering the full domain. 

For Stage 1 the domain is large but with small gradients in stress, therefore a mesh with large elements is used to reduce computational cost. Since Stage 2 will involve high stress gradients around the anchor, a mesh with significantly smaller elements is used, as shown in Figure \ref{fig:partdomain stage 2}. Therefore although the meshed domain in Stage 2 is smaller, the number of DoF of the background mesh is likely to be larger than that of Stage 1, hence a load step in Stage 2 is likely to be more expensive than Stage 1. The boundary conditions in Stages 1 and 2 are also different, for Stage 1 all surfaces, except the top surface, have roller boundary conditions applied so that the material can move freely under gravity. For Stage 2 all the boundary conditions, except the top surface, are fixed. The motivation for setting the boundary conditions as fixed is to prevent step changes in the material displacement for material that is either entering, or leaving, the partitioned domain. These boundary conditions can obviously be adjusted to suit other physical problems.  

\section{Modelling approach}\label{sec: Numerical simulations}

In this section we set out specific details of the model of a drag anchor pull and  validate the partitioned domain approach enabling long anchor pulls.  This is necessary in order to define the partitioned domain size for the full CPT-based predictive process to the prediction of anchor trajectories which is covered in Section~\ref{sec: CPT calibration}.



\subsection{Solution procedure}\label{sec: Anchor solution and recording procedure}

The procedure for predicting and recording the anchor trajectory is not trivial. This is particularly important in the context of validating against the centrifuge anchor pulls of Sharif ~\emph{et al.} \cite{Yaseen2024}, Section \ref{sec: anchor prediction}, as it is necessary to make sure the numerical procedure and recording of the anchor depth is as close to the experiment as possible.

\begin{figure}[ht!]
    \centering
    \begin{subfigure}[t]{1.0\textwidth}
        \centering
       \includegraphics[width=1\textwidth]{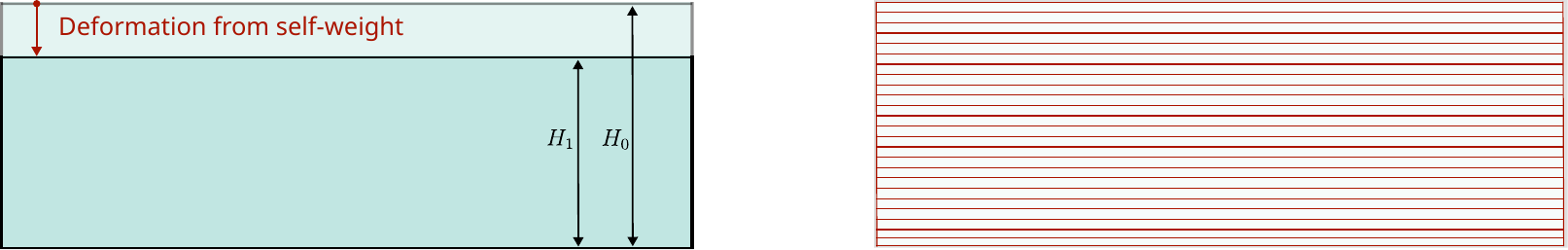}
        \caption{} \label{fig: non linear pro stage 1}
    \end{subfigure}%
   
    \begin{subfigure}[t]{1.0\textwidth}
        \centering
      \includegraphics[width=1\textwidth]{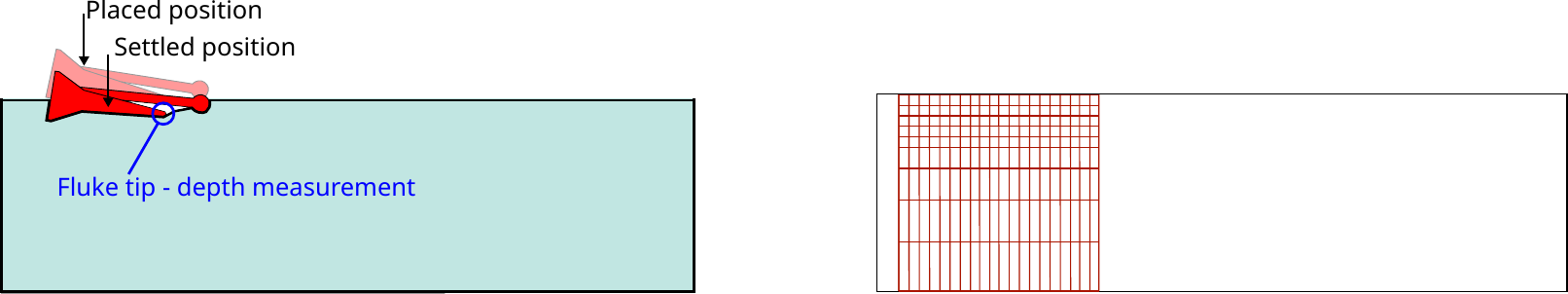}
        \caption{} \label{fig: non linear pro stage 2}
    \end{subfigure}

    \begin{subfigure}[t]{1.0\textwidth}
        \centering
      \includegraphics[width=1\textwidth]{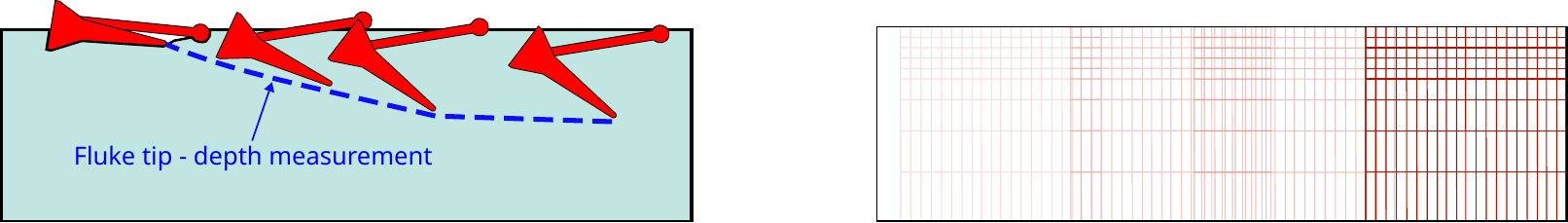}
        \caption{} \label{fig: non linear pro stage 3}
    \end{subfigure}
    \caption{Anchor solution and recording procedure:  The three stages and their corresponding meshes with: (a) Settling of the material under self-weight; (b) Placement and settling of the anchor on the surface; and (c) Anchor pull and recording of fluke tip position.}
\end{figure}

There are three key stages in an analysis: (a) initial stress state of the soil under self weight, (b) placing of the anchor on the soil surface and initial settlement under self weight and (c) anchor drag.  In addition, an analysis requires the generation of an initial distribution of material points, which could be uniform or graded depending on the problem being solved. 
The first stage, Stage~(a) shown in Figure \ref{fig: non linear pro stage 1}, 
is solved under quasi-static conditions to determine the initial stress state and displacement of the material points prior to an anchor test.
Stage (a) corresponds to Stage 1 of the partitioned domain approach, see Figure \ref{fig:partdomain stage 1}. Roller boundary conditions are on all surfaces, except the top. Since the material is homogeneous and has a flat upper surface, a stress variation only exists in the vertical direction. Therefore the corresponding mesh, shown on the right of Figure \ref{fig: non linear pro stage 1}, has element side lengths $0.1$ m in the vertical direction whilst being only $1$ element thick in the remaining directions.

The next two steps are shown respectively in Figures \ref{fig: non linear pro stage 2} and \ref{fig: non linear pro stage 3}. Both stages correspond to Stage 2 of the partitioned domain approach, see Figure \ref{fig:partdomain stage 2}. Fixed boundary conditions are applied on all surfaces, except the top where the material has a homogeneous Neumann interface and on the centreline/plane of symmetry where a symmetric (or roller) boundary condition is applied. On the right of Figures \ref{fig: non linear pro stage 2} and \ref{fig: non linear pro stage 3} an example of how the mesh through the domain as the anchor travels is also provided. The nature of the mesh used during (b) and (c) will depend on the anchor position and orientation. 

Stages (b) and (c) are solved using dynamic equilibrium, with the aim to represent as closely as possible the anchor pull procedure in the centrifuge \cite{Yaseen2024}. In the centrifuge the anchor is placed on the sand sample, then the rotational speed of the centrifuge is slowly increased to the required $g$ level (the spin up phase). Once at the correct $g$, the anchor is then pulled along the sample. The numerical Stage~(b) is designed to replicate the penetration of the anchor from self-weight during the spin up phase and is solved as a dynamic problem, which has the benefit of using the anchor's inertia to make the system determinate. The solution has the potential to be otherwise if solved as a quasi-static problem as not all parts of the anchor will be in contact with the sand from the first iteration of the first load step. These components will be unable to resist gravity, resulting in an ill conditioned system and large displacements making the non-linear solution difficult to solve. However, an issue with solving the problem dynamically is that there will be time dependent settling and vertical oscillation. To mitigate the effect of these oscillations on the drag embedment Stage~(c), the anchor is only pulled once all components of the anchor velocity are below $10^{-3}$m/s. In Stage~(c) the anchor is pulled via a linear truss element attached to the shank pad eye. This truss element is designed to replicate the steel wire that is attached to the shank in the centrifuge experiment. There are two key points here in terms of recording the kinematics and depth of the anchor penetration, which are based on the assumption that the fluke tip will achieve the greatest penetration depth:
\begin{enumerate}
    \item The zero offset for the anchor penetration depth is the depth of the fluke at the end of Stage~(b).
    \item The recorded depth of the anchor is \textit{always} the fluke tip.
\end{enumerate}
The purpose is to be consistent with the results of Sharif \emph{et al.} \cite{Yaseen2024}, where the depth is only recorded once the spin up phase is complete and the penetration depth is determined from the fluke tip.

\subsection{Constitutive behaviour}\label{sec: Constitutive behaviour}

The material model used to represent the sand in all of the analyses in this paper is a Hencky material (as described earlier) combined with a perfectly-plastic non-associated flow yield surface defined by a frictional cone with a Willam-Warnke \cite{WillamWarnke1974} deviatoric section\footnote{The Willam-Warnke deviatoric section \cite{WillamWarnke1974}, although originally designed for concrete, is very similar to other functions used for frictional soils, such as Lade-Duncan, Matsuoka-Nakai, Gudehus, Bhowmik-Long, etc. (see Bardet~\cite{Bardet1990} for a comparison). Critically it includes dependency on both the Lode angle and the intermediate principal stress.}. The constitutive parameters for this material model are obtained from Brinkgreve \emph{et al.}'s \cite{Brinkgreve2010} empirical relationships based on the relative density of the sand. The approach for going from a CPT trace to prediction of anchor penetration depth can be summarised as:
\begin{enumerate}
    \item Estimate a relative density for the CPT trace and obtain material properties from Brinkgreve \emph{et al.} \cite{Brinkgreve2010}.
    \item Run numerical CPT and obtain cone resistance trace.
    \item If the trace compares well, material properties are found; if the difference is large, increase or decrease the assumed relative density, update the material properties and go back to (2) to rerun the CPT.
    \item Use the calibrated material properties in the anchor simulations and determine the anchor trajectory.
\end{enumerate}

The scope of these validations is limited to sands with a homogeneous relative density and although the Brinkgreve \emph{et al.} \cite{Brinkgreve2010} empirical equations are used to provide the material data to obtain a CPT trace, it is acknowledged that a different set of values for the material data could be used to obtain a very similar trace. For the remaining sections all the material properties are defined with Brinkgreve's empirical equations, and hence are presented here for completeness (see Table~\ref{tab:Brinkgreve material properties}). The Willam-Warnke \cite{WillamWarnke1974} Lode angle dependent deviatoric section requires the specification of the deviatoric radius under triaxial extension normalised by that under triaxial extension, $\bar{\rho}_e$. In order to not introduce an additional constitutive parameter, in this paper $\bar{\rho}_e$ is obtained from the friction angle, $\phi$, by assuming that the Willam-Warnke section matches the Mohr-Coulomb criteria at the compression and extension meridians \cite{Coombs2011}
\begin{equation}
    \bar{\rho}_e = \frac{3-\sin(\phi)}{3+\sin(\phi)}.
\end{equation}
It is also assumed that the Young's modulus of the material varies with depth, where the initial Young's modulus for the point $p$ at initial depth $d_p$ is calculated using the formula by \cite{schanz2019hardening} assuming a negligible cohesion
\begin{equation}\label{eq: E}
    E_{50} = {E}_{50}^{r}\left(\frac{{{\sigma}_{v} {K}_{0}}}{{{p}^{r}}} \right)^{m_E}
    \qquad\text{with}\qquad 
    \sigma_v = g \rho d_p,
\end{equation}
where $K_0=1-\sin(\phi)$ is the coefficient of earth pressure at rest \citep{jaky1944coefficient}, $p^{r}=100$ kPa is the reference pressure, ${E}_{50}^{r}$ is the reference Young's modulus (see Table~\ref{tab:Brinkgreve material properties}) and $\sigma_v=g \rho d_p$ is the assumed vertical stress at $d_p$. Additional material parameters, which are constant for all relative densities are a Poisson's ratio of $\nu=0.25$ and a cohesion of $c =300$ Pa. As with Bird \emph{et al.} \cite{Bird2024}, the Young's modulus, and all other material properties, remain unchanged throughout the simulation. Since the model is linear-elastic, non-linear elastic response approximated by adopting the $E_{50}$ value as the constitutive model's Young's modulus, an approach also recommended in the PLAXIS handbook \cite{PLAXISV9}. 

\begin{table}[]
\centering
\caption{Empirical equations to determine the material data for relative density $R_D$, where the $R_D$ value in the equations varies between $0$ and $1$, i.e. a relative $82\%$ would be $R_D=0.82$. Gravitational acceleration is denoted by $g=9.81$ m/s$^2$. }
\begin{tabular}{rcll} \hline
    Reference Young's modulus ($E_{50}^{r}$) &$=$& $60R_D$                 & MPa\\
    Unit weight ($\rho^0\times g$)               &$=$& $ (19 + 1.6 R_D)-9.81$                     & kN/m$^3$\\
    Friction angle ($\phi$)          &$=$& $28 + 12.5R_D$      & $^{\circ}$\\
    Dilation angle ($\psi$)          &$=$& $\max(-2 + 12.5R_D,~0)$      & $^{\circ}$\\
    Stiffness exponent ($m_E$)        &$=$& $0.7 -0.3125 R_D$       & -\\ \hline 
\end{tabular}
 \label{tab:Brinkgreve material properties}
\end{table}

Friction between the rigid and the deformable bodies is modelled using Coulomb's friction law which captures the stick-slip behaviour in the tangential direction, described with the coefficient of friction $\mu$. The formulation is subject to the Karush-Kuhn-Tucker conditions which are only approximately enforced using a penalty formulation. The normal contact is also captured by a penalty formulation and therefore the Signorini-Hertz-Moreau conditions are also only approximated. The contact formulation, and implementation description, for contact between the GIMPs and the rigid body are described in 2D and 3D in \cite{Bird2024,bird2024dynamic}. The normal and tangential penalty constants for all simulations take the form,
\begin{equation}
    \epsilon_N = 50E\times A\quad\text{and}\quad \epsilon_T = 25E\times A,
\end{equation}
respectively, where the constants $50$ and $25$ have units $1$/length, in 2D $A$ is a characteristic length of the GIMP and in 3D $A$ is a characteristic area. It is important to note that other contact formulations exist, all with their own advantages and disadvantages, for example see the book by Wriggers \cite{wriggers2006computational}. However it was found in \cite{Bird2024,bird2024dynamic} that the penalty method, with these parameters, was a good compromise between accuracy, algorithmic stability and implementation difficulty. 

\subsection{Partitioned domain validation}\label{sec:partitioned domain validation}

The aim of this section is to validate the partitioned domain approach for drag anchors and to identify a partitioned domain size which is sufficiently large, for all relatively densities, so the boundaries do not influence the result. Therefore the drag anchor trajectory is recorded and compared for the different domain sizes and relative densities. Additionally the computational cost with drag distance is also recorded, although the primary aim of the method is to make the costs per load step invariant to the total length, it will be shown that a small increase in performance can be achieved by having a smaller partitioned domain\footnote{All timing tests were performed on 40 cores of a 2$\times$ AMD EPYC 7702 chipset using 200GB of 256GB RAM.}.

\subsubsection{Partitioned domain dimensions}

To validate the partitioned domain approach the dimensions for the partitioned domain need to be defined based on the geometry and location of the rigid body.

\begin{figure}[ht!]
    \centering
    \includegraphics[width=1.0\textwidth]{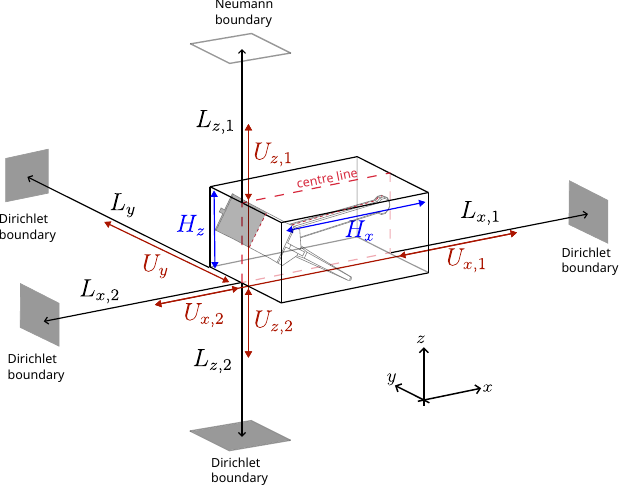}
    \caption{Dimensions for the partitioned domain.}
    \label{fig:partitioned-domain-dimensions}
\end{figure}

Figure \ref{fig:partitioned-domain-dimensions} illustrates the methodology behind defining the partitioned domain's lengths, and also the regions of uniform refinement, of the background mesh, around the partitioned domain. 

All the dimensions of the domain are defined in terms of the dimensions of a box which bounds the rigid body, as an example in Figure \ref{fig:partitioned-domain-dimensions} an AC-14 anchor has been used.
The steps to defining the domain are as follows:
\begin{enumerate}
    \item Identify the bounding box which surrounds the rigid body and identify its length $H_x$, width $H_y$ and height $H_z$.
    \item From the bounding box side lengths, $H_{(\cdot)}$, determine the distance from the edge of the box to the partitioned domain's boundary $L_{(\cdot),)\cdot{})}$.
    \item From $H_{(\cdot)}$ determine the regions of uniform refinement $U_{(\cdot),(\cdot)}$.
\end{enumerate}

For the anchor simulations in this paper the partitioned domain approach only affects the $x$ and $z$ directions with the $y$-dimensions of the mesh already defined in Section \ref{sec: Anchor solution and recording procedure} ($L_y = 5$ m and $U_y = 2$ m). From a numerical perspective the anchor problem is symmetric in the $x-z$ plane and no roll can occur. Therefore on the box the centre line of the anchor is marked, this corresponds to the plane where symmetric (roller) boundary conditions are applied. The distance of a partitioned domain boundary away from a face of the box, $L_{(\cdot),(\cdot)}$, is a function of the respective box size $H_{(\cdot)}$. Inside the box the mesh is uniform, additionally in the regions marked $U_{(\cdot),(\cdot)}$ (their size is also a function of $H_{(\cdot)}$) the mesh is also uniform.

For partitioned domain sizes the function which describes the region of uniform element sizes is kept constant
\[
U_{x,1} = U_{x,2} = \frac{1}{2}H_x\quad\text{and}\quad U_{z,1} = U_{z,2} = \frac{1}{2}H_z.
\]
Furthermore the domain size behind the anchor is also kept constant as this has minimal loading, $L_{x,2} = H_x/2$. The region above the anchor, $L_{z,1}$, is set to be sufficiently high to always include all the material above the anchor, and below the anchor $L_{z,2}$ always extends to the bottom of the full domain. Therefore the partitioned domain size which is varied is $L_{x,1}$,  for a loose sand ($32\%$ relative density) and a dense sand ($88\%$ relative density) with their setup and definition were outlined in Section \ref{sec: Constitutive behaviour}. The range of simulations is shown in Table \ref{tab:Anchor penetration: Partitioned domain lengths}.
\begin{table}[ht!]
    \centering
    \caption{Validation of partitioned domain: Partitioned domain lengths.}
    \begin{tabular}{r|c|c|c}
                 & test 1     & test 2     & test 3  \\ \hline
        $L_{x,1}$ (m) & $1.25H_x$ &  $1.5H_x$  & $2H_x$         
        \end{tabular}
    \label{tab:Anchor penetration: Partitioned domain lengths}
\end{table}

\subsubsection{Anchor design}\label{sec: Anchor design}
\begin{figure}[ht!]
    \centering
    \begin{subfigure}[t]{0.49\textwidth}
        \centering
       \includegraphics[width=1\textwidth]{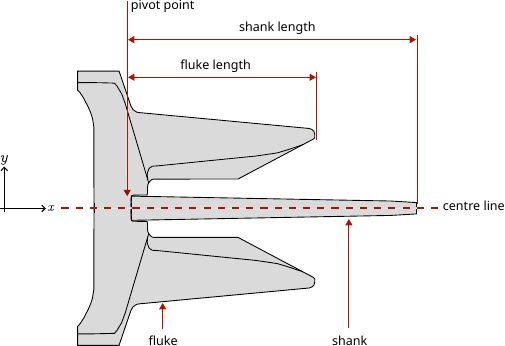}
        \caption{} \label{fig: anchor_schematic1}
    \end{subfigure}%
    \begin{subfigure}[t]{0.49\textwidth}
        \centering
      \includegraphics[width=1\textwidth]{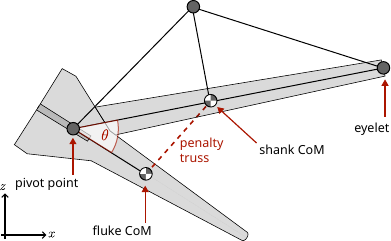}
        \caption{} \label{fig: anchor_schematic2}
    \end{subfigure}

    \begin{subfigure}[t]{1\textwidth}
        \centering
      \includegraphics[width=1\textwidth]{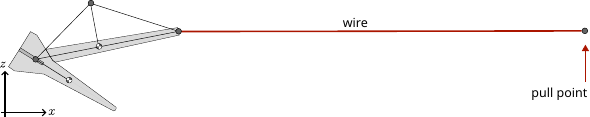}
        \caption{} \label{fig: anchor_schematic3}
    \end{subfigure}
    \caption{Anchor design: Top down schematic of the AC-14 anchor used for analysis is shown in (a) with the side profile in (b) and extended profile in (c).}
    \label{Fig:AC14 anchor schematic}
\end{figure}

The anchor centrifuge tests of Sharif \emph{et al.} \cite{Yaseen2024} use a model scale AC-14 anchor, here a full scale 8.7 tonne anchor is used, with a schematic for the anchor shown in Figures \ref{fig: anchor_schematic1} and \ref{fig: anchor_schematic2}.

The full anchor is represented\footnote{The anchor's fluke and shank are represented by a triangulated surface mesh with 372 and 188 facets respectively. The anchor STL is available as part of the paper's supplementary material.} but since a line of symmetry exists down the centre line of the anchor, Figure \ref{fig: anchor_schematic1}, it is only necessary to model half the material domain, achieved with a symmetric boundary condition. This assumption necessitates restraining the motion of the rigid body so it is only free to displace in the $x$ and $z$ (drag and penetration) directions and rotate about the $y$ (pitching) axis. To be consistent as only half the anchor is modelled, the anchor mass and rotational inertia also need to be halved. The full mass, rotational inertia and the position of the centre of mass (CoM) from the hinge point are provided in Table \ref{tab: anchor mass, intertia and length properties} with the anchor having a maximum opening angle of $\theta = 35^\circ$.

\begin{table}[ht!]
\centering
\caption{Anchor design: Total mass, rotational inertia (rot. inertia) and geometric properties of the anchor components, where the CoM position and length and measured from the pivot point.}
\label{tab: anchor mass, intertia and length properties}
\begin{tabular}{r|c|c|c|c}
      & mass (kg)  & rot. inertia (kgm$^2$) & CoM position (m) & length (m) \\ \hline
fluke & 6583.2 & 1100               & 0.131                   &  1.7 \\ \hline
shank & 2116.7 & 1350               & 1.272                   & 3.3 \\ 
\end{tabular}
\end{table}

The wire, shown by the thick red line in Figure \ref{fig: anchor_schematic3}, is modelled with a truss element and has a length of $9.5$ m, stiffness of $10^{10}$ N/m. The end of the truss element not attached to the anchor, the pull point, has a prescribed velocity of $\{v\}=\{0.1 \quad 0 \quad 0\}^T~\text{m/s}$. The initial height of the pull point is set to coincide with the initial surface of the soil ($z=8$m in all of the examples in this paper). Last, the coefficient of friction acting at the sand-anchor interface is $0.3$.


\subsubsection{Material point discretisation}

To validate the partitioned domain approach and determine an appropriate domain size, two different relative density sands were considered, $32\%$ and $82\%$, with the material properties defined by the equations in Section~\ref{sec: Constitutive behaviour}. The mass and stiffness ghost stabilisation \cite{coombs2022ghost} parameters are set to $\gamma_M=\rho/4$ and $\gamma_K=E/10$, respectively, where $\rho$ and $E$ are the volume weighted average density and Young's modulus of the material points that occupy the elements that share the element boundary where the stabilisation is applied.  
The stages to solving and recording the anchor trajectory defined in Section \ref{sec: Anchor solution and recording procedure}. In total the anchor was dragged a distance of $7$m.

\begin{figure}[ht!]
    \centering
    \includegraphics[width=0.6\textwidth]{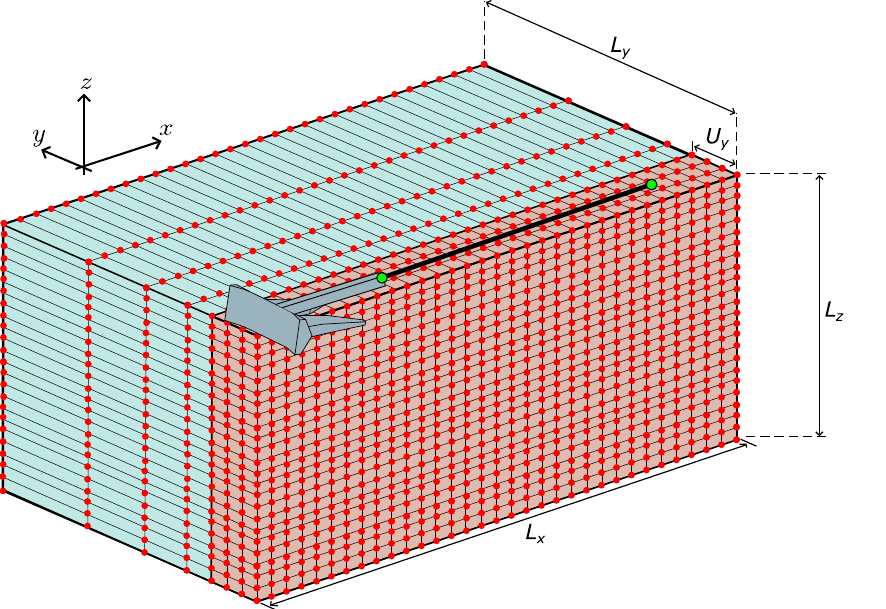}
    \caption{Partitioned domain validation: Initial distribution of the background grid nodes for the GIMP setup.}
    \label{fig:anchor penetration - initial GIMP setup}
\end{figure}

As discussed in Section \ref{sec: Partitioned domain approach}, an initial distribution of GIMPs for the full domain is required. For this validation, a Cartesian mesh is generated which is sufficiently large to encompass the full anchor trajectory with GIMPs which are evenly distributed in a $2\times2\times2$ grid within each element. An indicative mesh is shown in Figure \ref{fig:anchor penetration - initial GIMP setup} which shows the exterior mesh. Figure \ref{fig:anchor penetration - initial GIMP setup} shows two regions: the red region contains uniform elements with side lengths of $0.1$ m, and in the blue region the element side lengths are increased by a factor of $1.3$, starting with side length $0.1$ m, in the $y$-direction between adjacent elements. The total domain size is $L_x = 25$ m,  $L_y = 5$ m, $L_z = 8$ m with $U_y = 2$ m. Figure~\ref{fig:anchor penetration - initial GIMP setup} provides a schematic of the problem setup prior to the anchor pull. 

\subsubsection{Results and discussion}

\begin{figure}[ht!]
    \centering
    \begin{subfigure}[b]{0.49\linewidth}
        \centering
        \includegraphics[width=\linewidth]{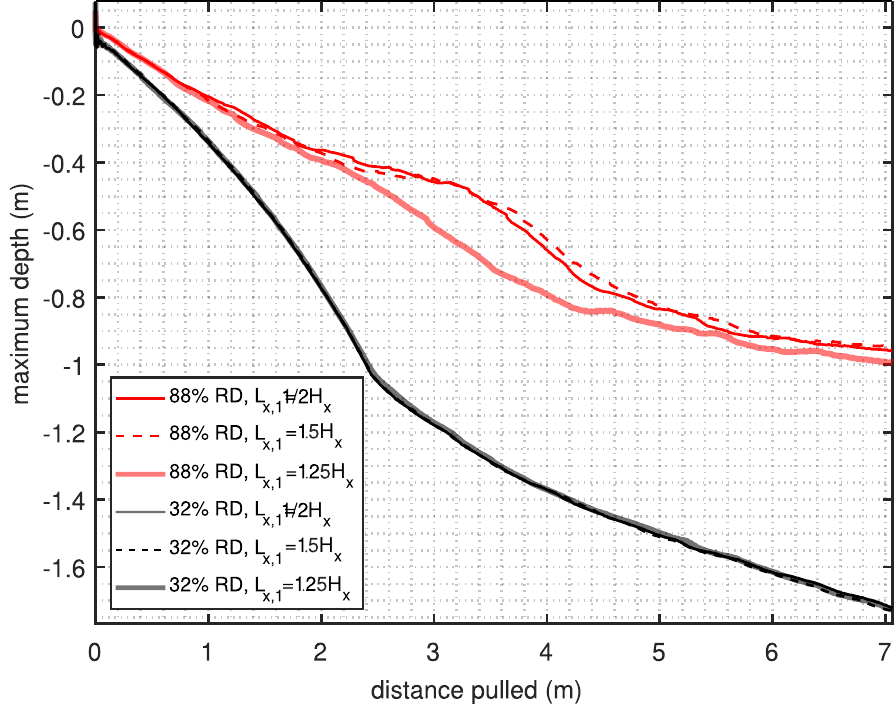}
        \caption{}
        \label{fig: pd depth test}
    \end{subfigure}
    \hfill
    \begin{subfigure}[b]{0.49\linewidth}
        \centering
        \includegraphics[width=\linewidth]{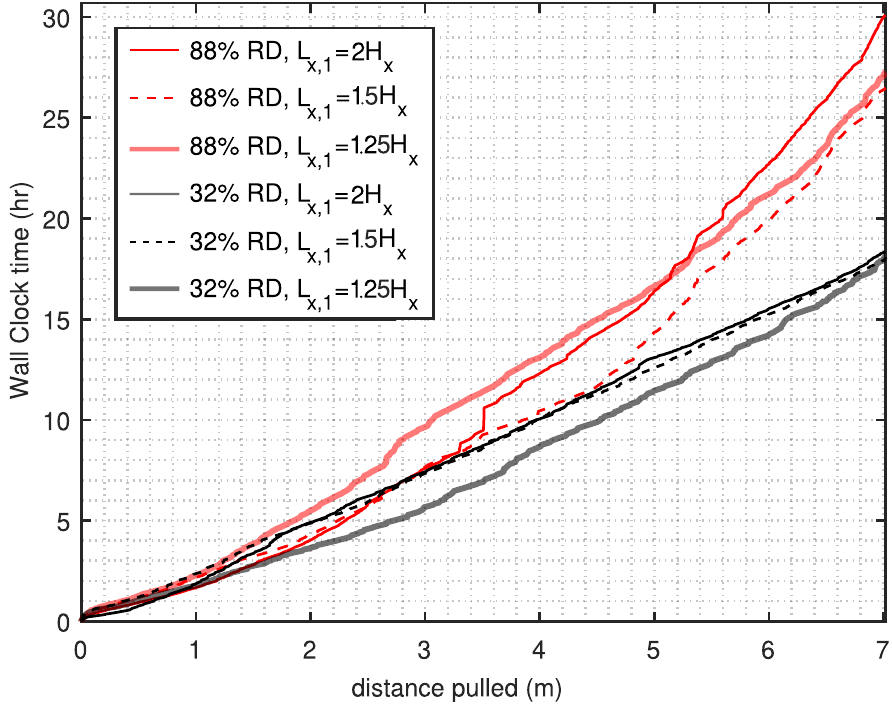}
        \caption{}
        \label{fig: pd timing test}
    \end{subfigure}
    \caption{Partitioned domain validation: The penetration depth of the anchor with drag distance is shown in (a) with the corresponding cumulative time for the simulation shown in (b).}
\end{figure}

The results for the three partitioned domain tests, Table \ref{tab:Anchor penetration: Partitioned domain lengths} are shown in Figures \ref{fig: pd depth test} and \ref{fig: pd timing test}. The first observation of the anchor trajectory, Figure \ref{fig: pd depth test}, is that the $32\%$ relative density sand has consistent results for all three domain sizes whereas $82\%$ agrees well for the two larger domain sizes with the smallest domain size an outlier. The results show that having $L_{x,1}=2H_x$ is sufficient for the range of materials considered here, with $L_{x,1}$ having a strong effect on the higher relative densities. For the $32\%$ relative density it is possible that a smaller partitioned domain size could be considered.
Ideally the reaction force, generated by the anchor's motion, at the boundary associated with $L_{x,1}$, Figure \ref{tab:Anchor penetration: Partitioned domain lengths}, should be zero. However, the probable cause of the greater dependence on $L_{x,1}$ for the higher relative density is its higher Young's modulus and the general higher strength of the material. Hence a larger $L_{x,1}$ is required. In summary Figure \ref{tab:Anchor penetration: Partitioned domain lengths} validates the domain approach, the results converge with increasing domain size for a range of relative densities.

The timing results are shown in Figure \ref{fig: pd timing test}, the primary observations are: (i) The $88\%$ relative density sand is consistently more expensive than the $32\%$; (ii) For this range of relative density there is not a correlation between the computational cost and the domain size. The reason for the lack of correlation is that the difficulties of solving a large deformation contact problem with a iterative solver dominate the solution time, rather than the number of GIMPs and Number of Degrees of Freedom (NDOF). This presents itself either as large number of iterations for the linear solver, a large number of NR-iterations, or, having to restart the load step with a reduced time increment. An example of such a case is at distance $3.5$ m of $88\%$ $L_{x,1} = 2H_x$ where the wall clock time continuous to increase but the distance pulled is unchanged. Combining these observations with the observations of Figure \ref{fig: pd depth test} shows that for the necessary domain size length of $L_{x,1}=2H_x$ there is no compromise on in terms of solution time. Therefore this domain size will be used for the remaining anchor simulations.

\begin{figure}[ht]
          \centering
          \includegraphics[width=0.8\linewidth]{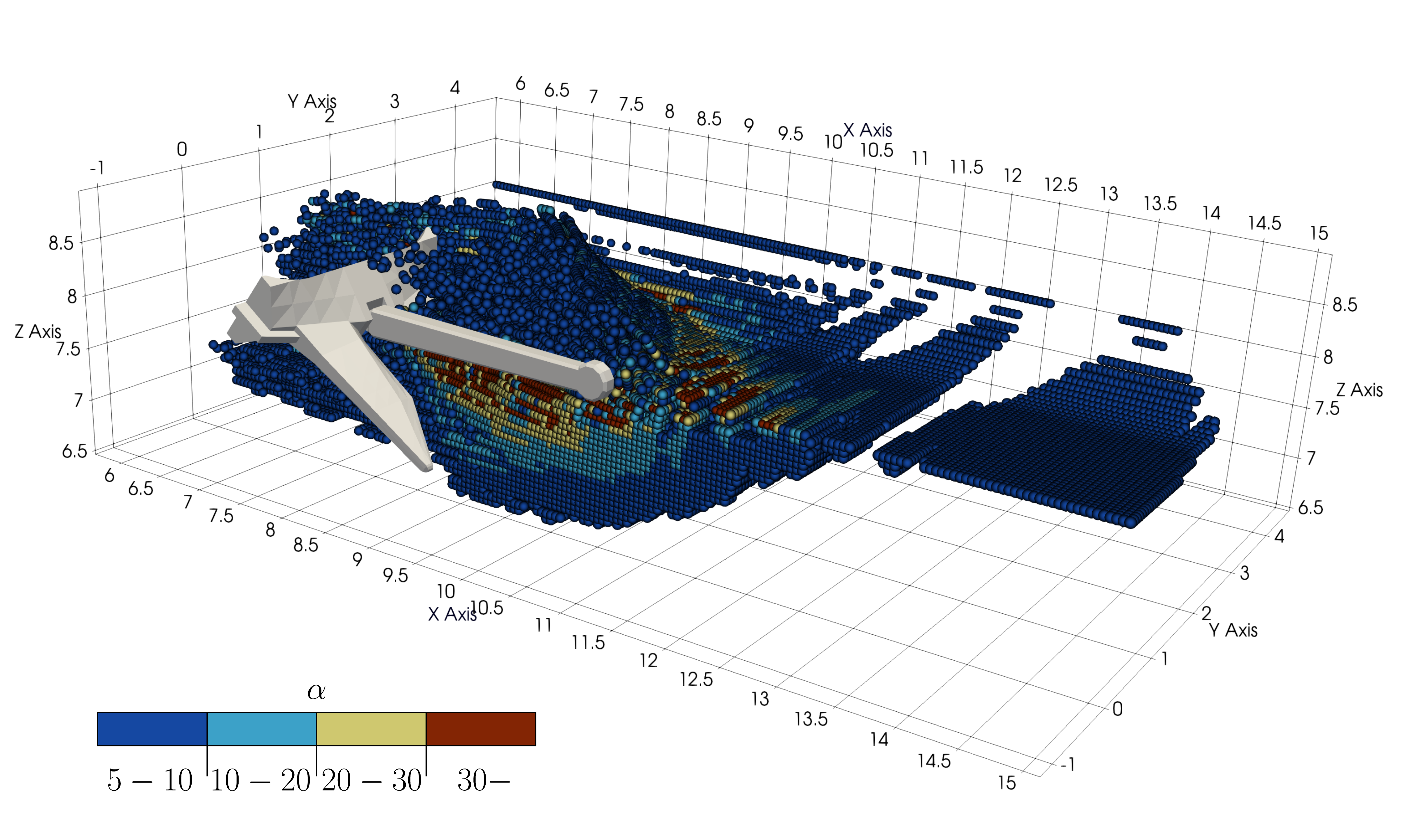}
          \caption{Validation of partitioned domain: Visualisation of the anchor at a pull distance of $7$ m when $L_{x,1}=2H_x$.}
          \label{fig:rd88_2_stress_Ld2}
\end{figure}

\begin{figure}[h!]
          \centering
          \includegraphics[width=0.8\linewidth]{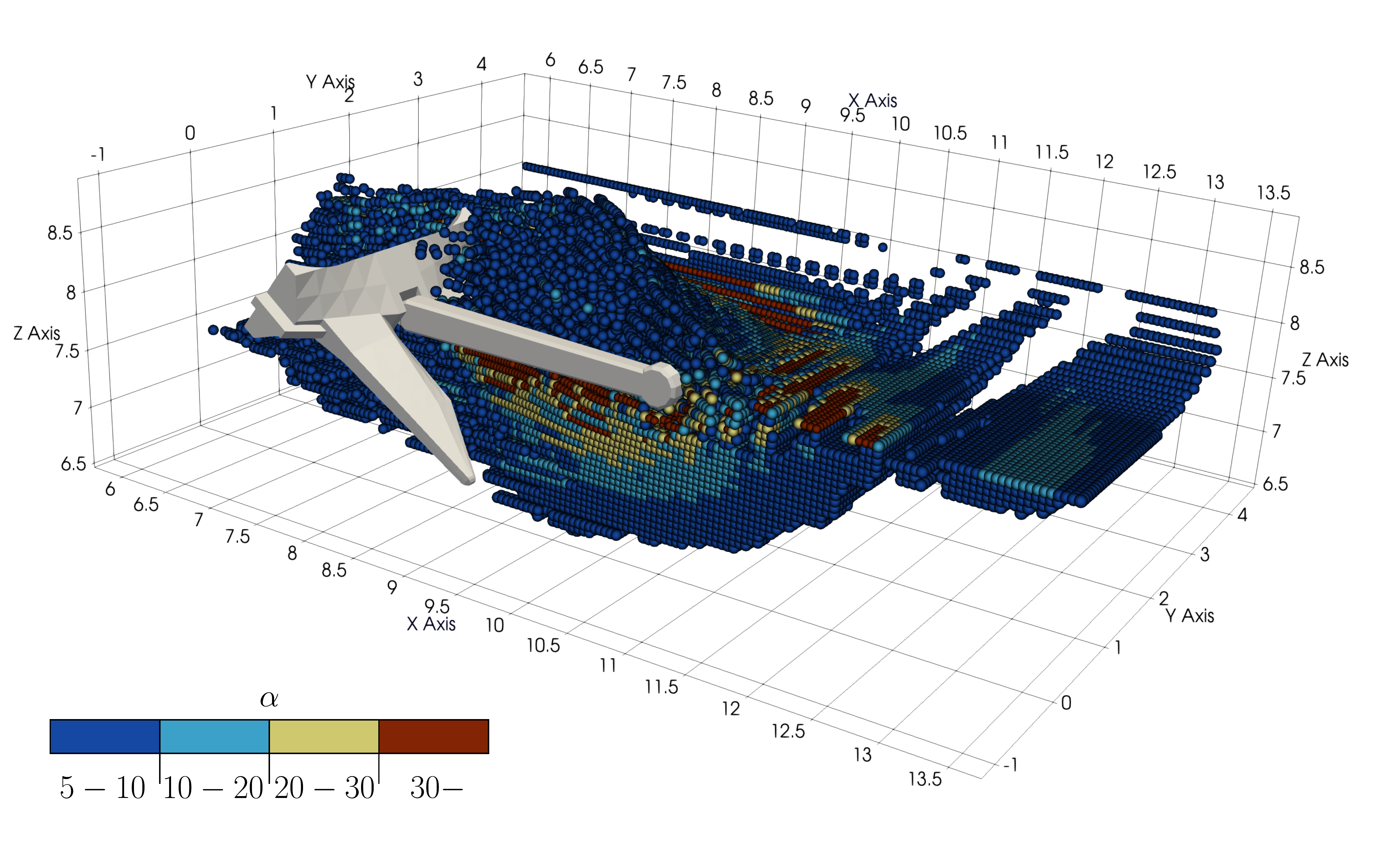}
          \caption{Validation of partitioned domain: Visualisation of the anchor at a pull distance of $7$ m when $L_{x,1}=1.5H_x$.}
          \label{fig:rd88_2_stress_Ld15}
\end{figure}

\begin{figure}[h!]
          \centering
          \includegraphics[width=0.8\linewidth]{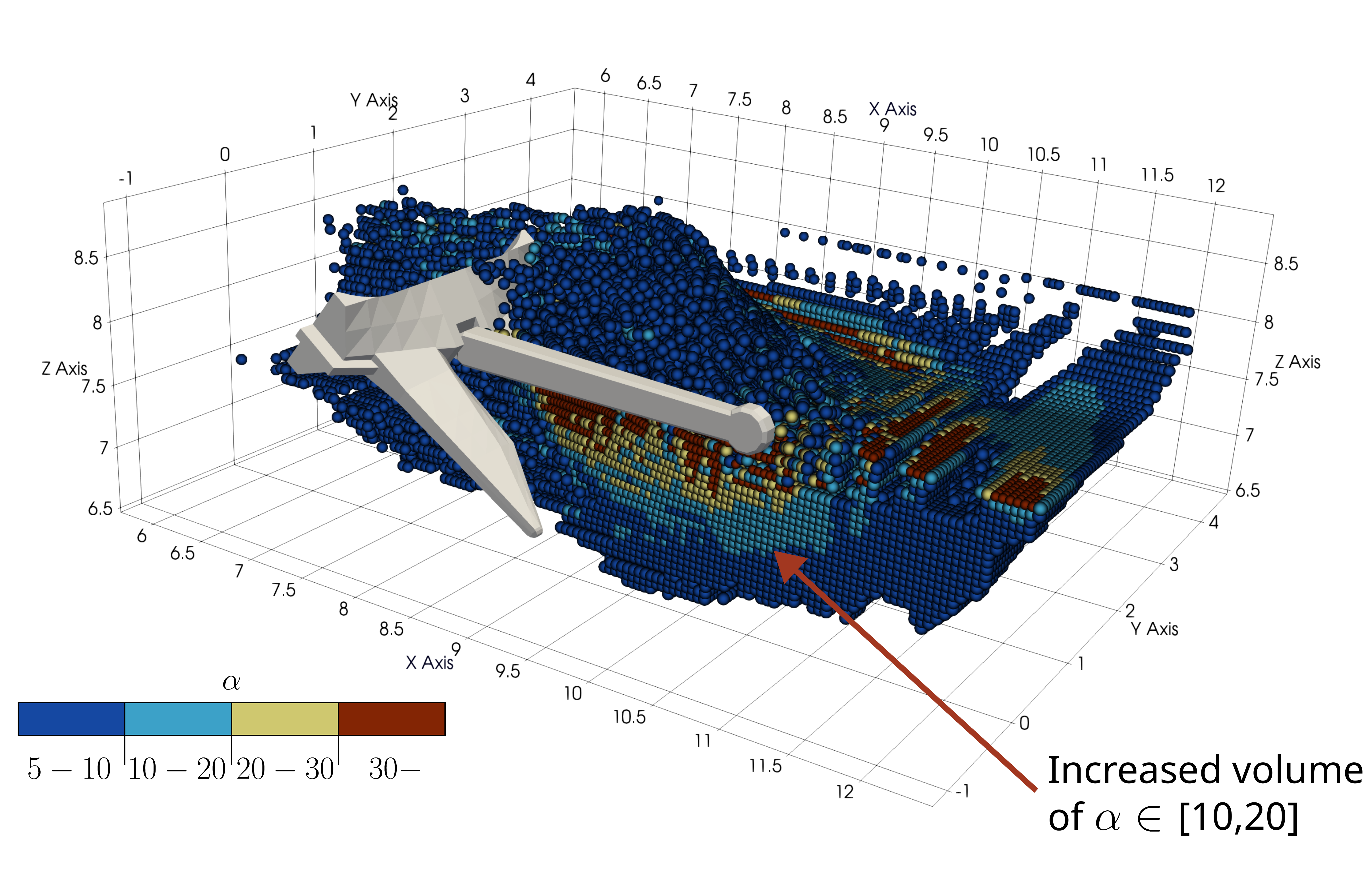}
          \caption{Validation of partitioned domain: Visualisation of the anchor at a pull distance of $7$ m when $L_{x,1}=1.25H_x$.}
          \label{fig:rd88_2_stress_Ld125}
\end{figure}

Inspection of the hydrostatic stress for the 88\% relative sand provides evidence as to why the depth profiles in Figure \ref{fig: pd depth test} vary with $L_{x,1}$ but also why increasing $L_{x,1}$ leads to convergence in the depth profile. To inspect the hydrostatic stress state the following parameter, similar to \cite{dao2024numerical}, is used
\begin{equation}
    \alpha = \frac{\Delta I_1}{I_{1,0}}
\end{equation}
where $I_{1,0}$ and $\Delta I_1$ are the initial and change in the first stress variant of Cauchy stress. Three stress plots are presented Figures \ref{fig:rd88_2_stress_Ld2}, \ref{fig:rd88_2_stress_Ld15}, and \ref{fig:rd88_2_stress_Ld125}, which respectively correspond to $L_{x,1}=2H_x$, $L_{x,1}=1.5H_x$ and $L_{x,1}=1.25H_x$. The plots show the material points and are coloured with respect to four $\alpha$ ranges, $\alpha\in[5,10]$, $\alpha\in[10,20]$, $\alpha\in[20,30]$ and $\alpha>30$.
The Figures show that as $L_{x,1}$ is increased the $\alpha$-plots become more similar.  $L_{x,1}=2H_x$ and $L_{x,1}=1.5H_x$ have the most similar $\alpha$-plot, with a difference of $0.5H_x$, whereas the comparison between $L_{x,1}=1.5H_x$ and $L_{x,1}=1.25H_x$ shows a significant difference in $\alpha$-plots despite a smaller change of $0.25H_x$
At $x=12$ m $L_{x,1}=1.25H_x$ a vertical front for the range $\alpha\in[5,10]$. Additionally high $\alpha$ values are observed along the surface in front of the anchor for $L_{x,1}=1.25H_x$ which are not observed for $L_{x,1}=1.5H_x$. Below the eyelet of the shank for $L_{x,1}=1.25H_x$ (marked with the red arrow) there also appears to be an increase in the volume of for $\alpha\in[5,10]$ compared to the other two cases. The similarities and differences in the surface for $\alpha>2$ for the three figures correspond to the similar behaviour in the depth profile in Figure \ref{fig: pd depth test}, with $L_{x,1}=2H_x$ and $L_{x,1}=1.5H_x$ being the most similar. Since the results for depth are still marginally changing from Figure \ref{fig:rd88_2_stress_Ld15} to Figure \ref{fig:rd88_2_stress_Ld2} the value of $L_{x,1} = 2H_x$ is used for the Anchor prediction validation in Section \ref{sec: anchor prediction}.

\section{CPT-based anchor penetration prediction}\label{sec: CPT calibration}

In this section CPT results are used to calibrate the material parameters in the numerical model by comparing the numerical trace for cone resistance with depth with the traces obtained by the centrifuge CPT for four sand densities,  \cite{davidson2022physical,Cerfontaine2020}, which were used in the centrifuge anchor pulls. The implicit-axisymmetric CPTs were validated in \cite{Bird2024,coombs2025aggregated} for a range of numerical parameters, therefore for this set of computations the numerical parameters are kept constant whilst the relative density is varied. Each numerical CPT therefore had a wall clock time of approximately 15 minutes \cite{Bird2024}, run as a serial computation in MATLAB\footnote{All computations were computed on a machine with Intel Core i7-8665U (4 cores, 8 threads) with 32 GB RAM running Ubuntu 24.04 LTS.}, meaning the material calibration was a relatively fast process.

The implicit quasi-static GIMPM analysis of a CPT by Bird~\emph{et al.} \cite{Bird2024} demonstrated good agreement between numerical predictions and experimentally measured $q$ values across of a range of relative densities where the relative density was known \emph{a priori}. However, in non-lab conditions the relative density of the sand is not known and often only a $q$-trace is provided (in reality offshore geotechnical data availability can be even more limited \cite{macdonald2023depth}). This necessitates a framework which first estimates the material properties from limited CPT data, and second, uses these estimated properties to predict anchor trajectory and ultimate penetration depth at the CPT location. 

The material properties are defined using Brinkgreve's \emph{et al.} \cite{Brinkgreve2010} empirical relations, see Section \ref{sec: Constitutive behaviour}. The relative density parameter is used to generate the material properties to obtained agreement between the numerical and experimental results, however the general trace of the CPT is not unique for a given set of parameters and a combination of other parameters could also generate similar results. Hence the relative density determined using the numerical CPT is not a determination of the true relative density, this is particularly important to note at the extreme values of relative density since it is likely in these regions the empirical relations of Brinkgreve's \emph{et al.} \cite{Brinkgreve2010} were dervied from fewer data points.

The mass and stiffness ghost stabilisation \cite{coombs2022ghost} parameters are set to $\gamma_M=\rho/4$ and $\gamma_K=E/10$, respectively, where $\rho$ and $E$ are the volume weighted average density and Young's modulus of the material points that occupy the elements that share the element boundary where the stabilisation is applied.

\subsection{CPT Calibration}\label{sec: CPT calibration}

\subsubsection*{Setup}

\begin{figure}[ht!]
    \centering
    \begin{subfigure}[b]{0.35\linewidth}
        \centering
        \includegraphics[width=\linewidth]{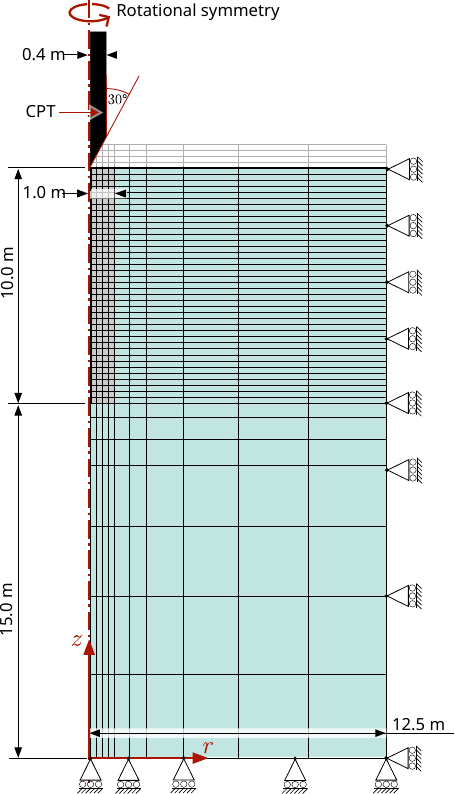}
        \caption{}
        \label{fig: CPT description}
    \end{subfigure}
    \hfill
    \begin{subfigure}[b]{0.64\linewidth}
        \centering
        \includegraphics[width=\linewidth]{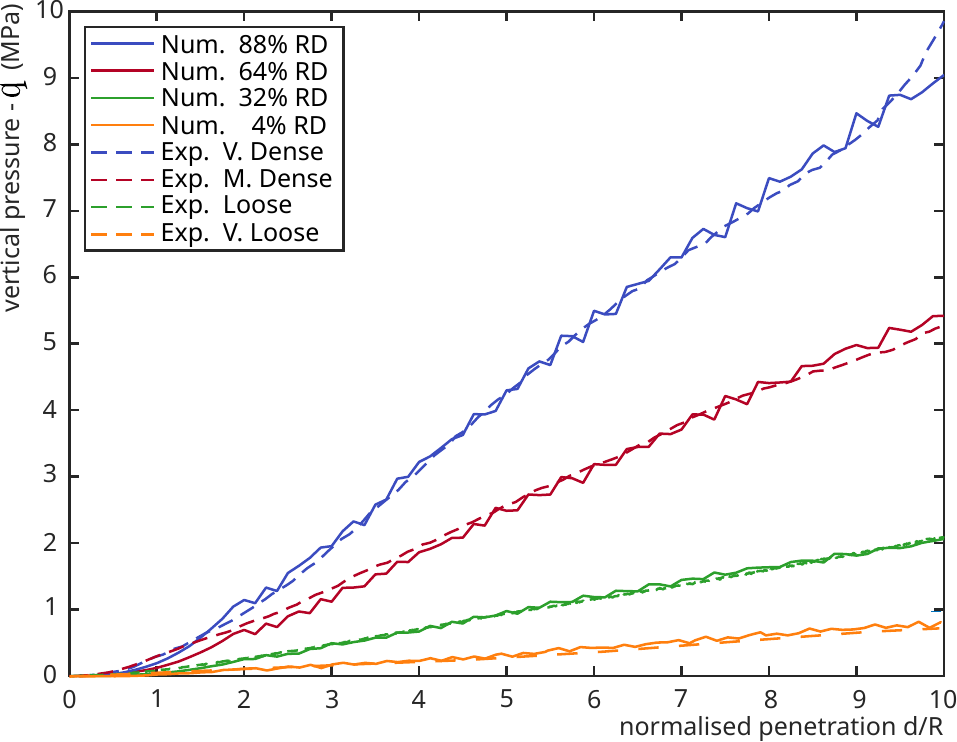}
        \caption{}
        \label{fig: CPT results}
    \end{subfigure}
    \caption{CPT calibration: (a) is a diagram of the overall geometry and setup of the numerical CPT with an indicative, and (b) showing the final calibration of numerical cone resistance, against normalised penetration, for a range of laboratory CPT results.}
\end{figure}

The material equations are described in Section~\ref{sec: Numerical simulations}, here the geometry, mesh, boundary conditions, and GIMP properties are defined for the CPT simulations. A diagram of the overall setup is provided in Figure \ref{fig: CPT description} which shows the geometry, boundary and indicative mesh. The CPT has a radius of $0.4$m and will penetrate into the sand $4.0$ m, the coefficient of friction between the CPT and sand is $0.3$. In Figure \ref{fig: CPT description} the region in which the CPT will penetrate is marked in grey, it has a refined and uniform mesh where the elements have a side length of $0.1$ m. Outside the grey region the element side length of adjacent element is increased by a factor of $1.3$. Initially each element throughout the mesh has $3^2$ material points. 

The simulation procedure is to first apply the gravitational load, then over $100$ load steps the CPT is moved vertically downwards over the distance of $4.0$ m. The load on the tip of CPT, from contact with the GIMPs, is recorded and then divided by the cross sectional area to give the cone tip resistance $q$ MPa. The cone resistance is then used for the calibration procedure outlined in Section \ref{sec: Numerical simulations} for the four sands used by Sharif \emph{et al.} \cite{Yaseen2024} in the anchor centrifuge tests; \textit{very dense}, \textit{medium dense}, \textit{loose} and \textit{very loose}.

\subsubsection*{Results and discussion}

The result of the calibration procedure for the four  sand Relative Densities (RDs) is provided, with the results shown in Figure \ref{fig: CPT results} with the results summarised in Table \ref{tab: exp. to. num material data}. Inspection of the Figure \ref{fig: CPT results} shows the relatively simple material model is able to sustain a very good level of agreement of a large range of normalised penetration, no adjustment of the material parameters is required over the penetration depths considered. This will also be validated for the anchor trajectories in Section \ref{sec: anchor prediction} for a range of sands considered here.

\begin{table}[ht!]
    \centering
    \caption{CPT calibration: Summary of the calibrated Relative Density (RD) values versus the experimental description.}
    {\small  \begin{tabular}{r|c|c|c|c}
        {Exp. description} & Very dense & Medium dense & Loose  & Very loose \\ \hline
        Calibrated RD & $88\%$     & $64\%$       & $32\%$ & $4\%$
    \end{tabular}}
        \label{tab: exp. to. num material data}
\end{table}

\subsection{Anchor penetration prediction}\label{sec: anchor prediction}

\subsubsection*{Setup and scope}

The material properties for each of the densities in Table \ref{tab: exp. to. num material data} are used to calculate the anchor trajectories and final penetration depths of the AC-14 anchor, completing the prediction procedure which started with a physical cone resistance trace. The coefficient of friction acting between the anchor and the sand is set to $0.3$.

The final penetration depths of the AC-14 are compared to those obtained in the centrifuge by Sharif ~\emph{et al.} \cite{Yaseen2024} for four sands of homogeneous relative density determined in the previous section.

\begin{figure}[ht!]
    \centering
    \includegraphics[width=0.6\textwidth]{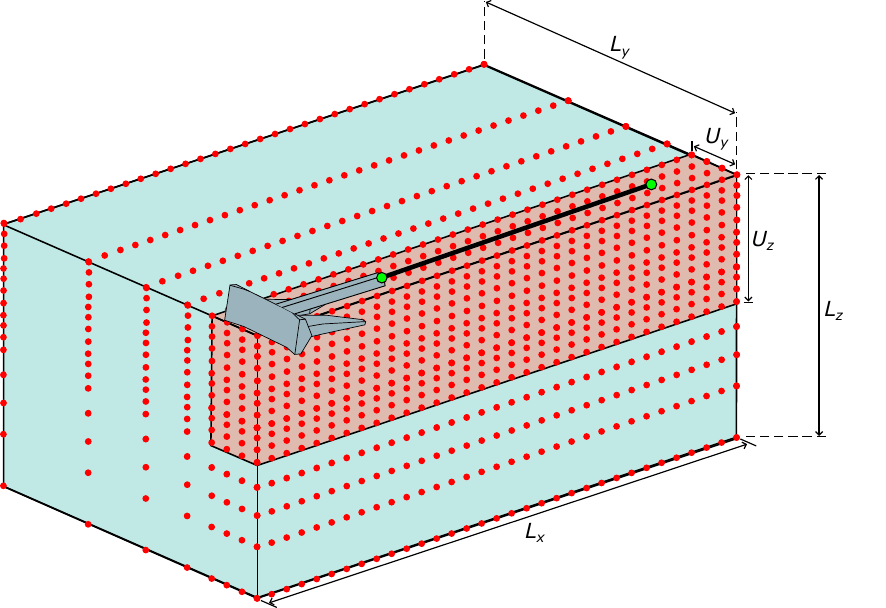}
    \caption{Anchor penetration: Schematic of the initial distribution of background grid nodes.}
    \label{fig:anchor penetration - initial GIMP setup with Uz}
\end{figure}

The setup for this problem is similar to the validation of the partitioned domain, Section \ref{sec:partitioned domain validation}, the only difference is that the total length of the domain is set to $L_x = 50$ m and the volume of the initially uniformly distributed GIMPs is reduced from being the total depth of the domain to $U_z = 4$ m, Figure \ref{fig:anchor penetration - initial GIMP setup with Uz}, such that $L_z = 8$ m. In the $z$-direction, below the uniform region, the adjacent node spacing was increased by a factor of $1.3$. 

\begin{table}[]
 \caption{Anchor penetration: Partition domain parameter values.}
    \centering
    {\small \begin{tabular}{l|c|c|c|c|c|c}
         Dimension &  $U_{x,1}$, $U_{x,2}$ & ~~$U_y$~~ & $U_{z,1}$, $U_{z,2}$ & $L_{x,1}$ & $L_{x,2}$ & ~~$L_{y}$~~\\ \hline 
         Value (m) & $H_x/2$ & $2$ & $H_x/2$ & $2 H_x$ & $0.5 H_x$ & $5$
    \end{tabular}}
    \label{tab: partitioned domain parameteres}
\end{table}

Given the results of the partitioned domain validation, the parameters for the partitioned domain and the volume of uniform refinement about the anchor are provided in Table \ref{tab: partitioned domain parameteres}. The material model, and setup, is consistent with the CPT calibration in Section \ref{sec: CPT calibration}, no alteration to the relative density, and therefore material parameters, are made from those summarised in Table \ref{tab: exp. to. num material data} and there are no further parameters or complexity added to the material model. Last, the simulations are stopped if over $2$m of drag the depth of the anchor does not increase by $10$mm. $2$m was considered as it is greater than the length of the fluke, $1.7$m, and so therefore represents the total process of sand moving over and around the fluke. The depth value is taken as this is considered to be small enough, with respect to $2$m, that no more significant penetration will occur.  

\subsubsection*{Results and discussion}

\begin{figure}
    \centering
    \includegraphics[width=\textwidth]{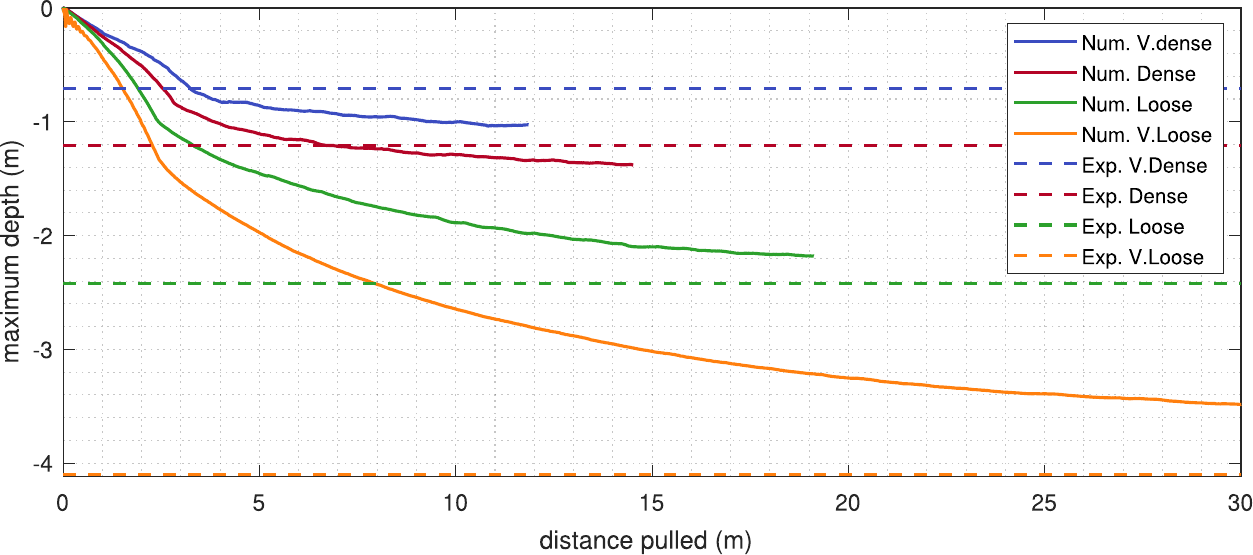}
    \caption{AC-14 anchor penetration: Comparison of anchor trajectories with the final depth of the experimental results.}
    \label{fig:Anchor penetration}
\end{figure}

The results for the anchor penetration are presented and compared to the final experimental penetration depths, are shown in Figure \ref{fig:Anchor penetration}. It must be emphasised that the results in this section are predictions of the anchor behaviour; no numerical or material data tuning has been performed to get results that agree. Furthermore the trajectory of the anchor is purely driven by the response of the sand, and the sand response is purely driven by the anchor.

In general the trend of the results matches that of the experimental data, for a high relative density the penetration is less and the difference between the final penetration depth increases with decreasing relative density. Furthermore the comparison of the results for the middle relative densities are very good, for the dense sand the absolute difference is $0.1$m and the loose sand is $0.2$m (both correspond to $\approx 8\%$ relative difference). Good agreement is also achieved for the very loose sample with a predicted final penetration of $3.5$ m compared to the experimental result of $4.1$ m, corresponding to an difference of $14\%$, or relative to the anchor length, a fraction of $0.35$ of the fluke length. Last, the $88\%$ has an absolute difference of $0.3$m, but with a higher relative difference of $43\%$, compared to the final penetration depth. The reason for this discrepancy can be explained by the differences in the kinematics between the experimental and numerical anchor for the very dense sand. 

\begin{figure}[ht!]
    \centering
    \includegraphics[width=1.0\linewidth]{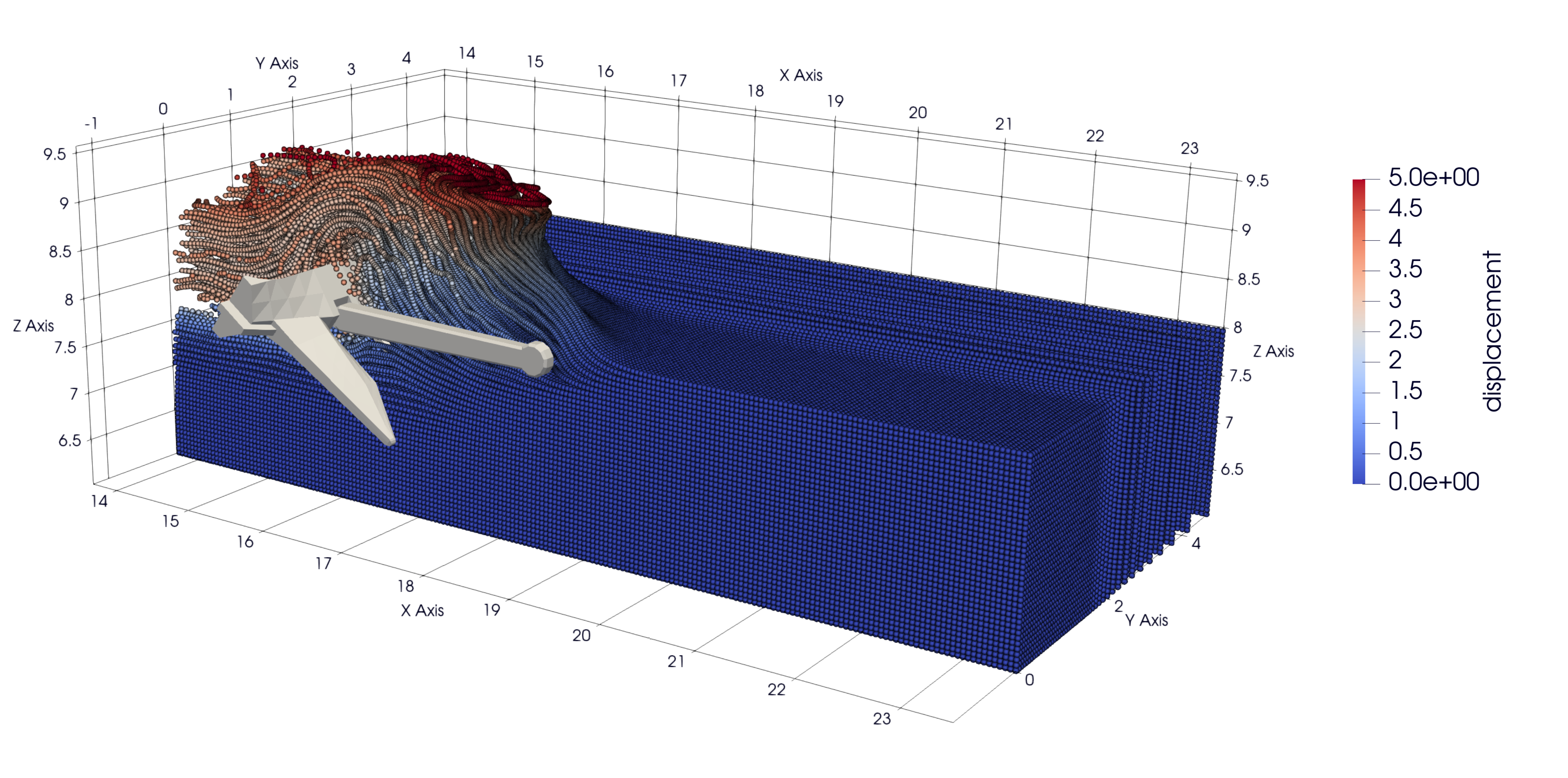}
    \caption{AC-14 anchor penetration: anchor position for the very dense sand at a $12$m pull distance.}
    \label{fig:rd_82_full_plot.pdf}
\end{figure}

Figure \ref{fig:rd_82_full_plot.pdf} shows the AC-14 anchor pulled in a sand of 88\% relative density at a 12 m pull distance, the shank is laying nearly horizontally at $0.5^\circ$ to the undeformed surface of the sample whilst the fluke is at $36.14^\circ$ degrees to the horizontal. The result is that nearly all the penetration is due to the flue inclination. This is in contrast to the experimental results where the penetration is due to burying of the fluke at a near horizontal angle of $5\%$ to the surface. The difference in the kinematics could potentially be caused by the assumption of a frictionless pin between the fluke and the shank in the numerical model. The physical tests are conducted in a geotechnical centrifuge which place the anchor assembly, and therefore the joint between the fluke and the shank, under large forces which will result in some frictional resistance to the anchor opening. There is also the potential for sand to be \emph{trapped} in the connection during the centrifuge test that again could provide additional resistance to shank-fluke rotation.

\begin{figure}[ht!]
    \centering
    \includegraphics[width=1.0\linewidth]{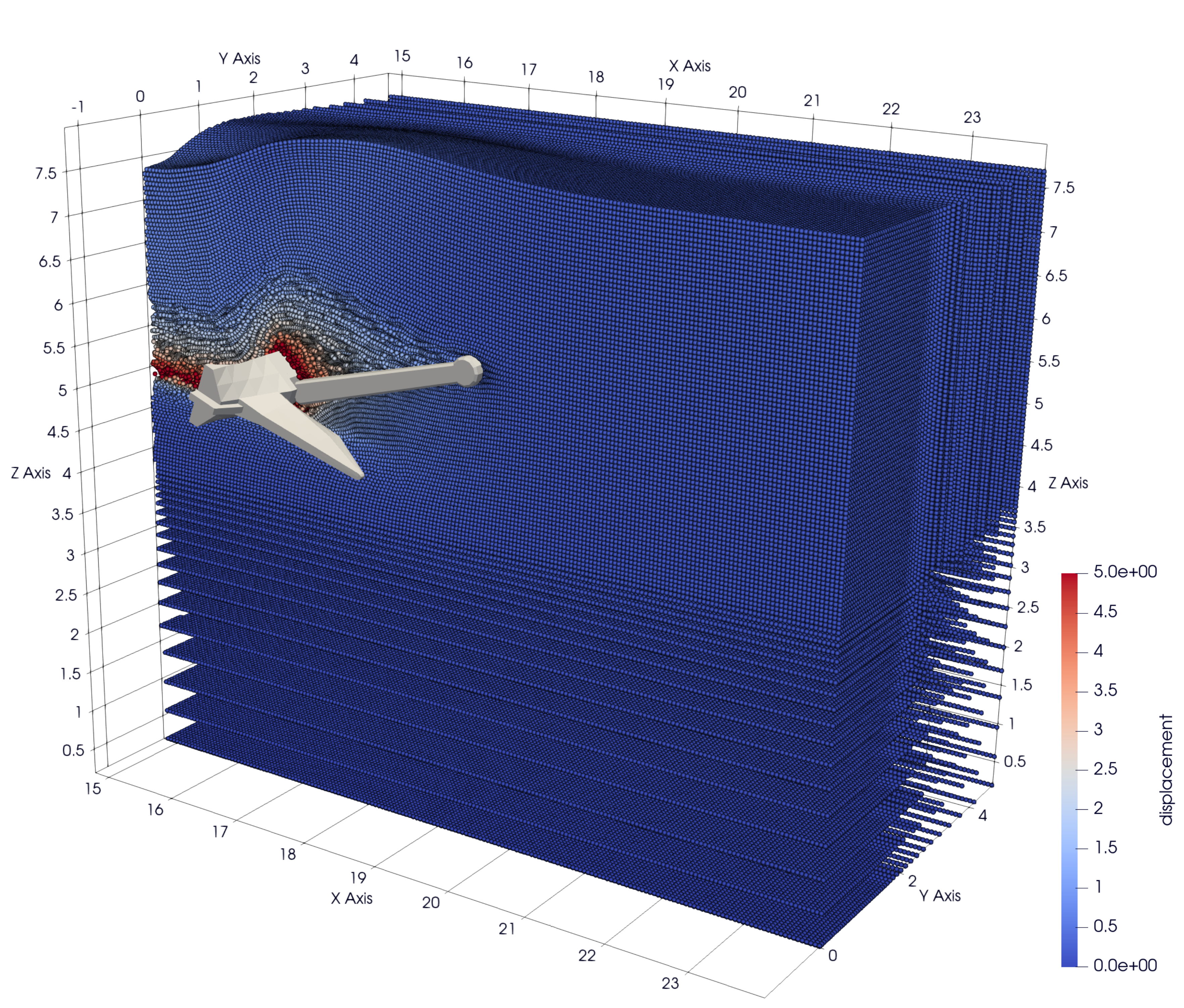}
    \caption{AC1-14 anchor penetration: anchor position for the very loose sand at a $17$m pull distance.}
    \label{fig:rd_4_full_plot.pdf}
\end{figure}

\begin{figure}
    \centering
    \begin{subfigure}[b]{0.8\linewidth}
        \centering
        \includegraphics[width=\linewidth]{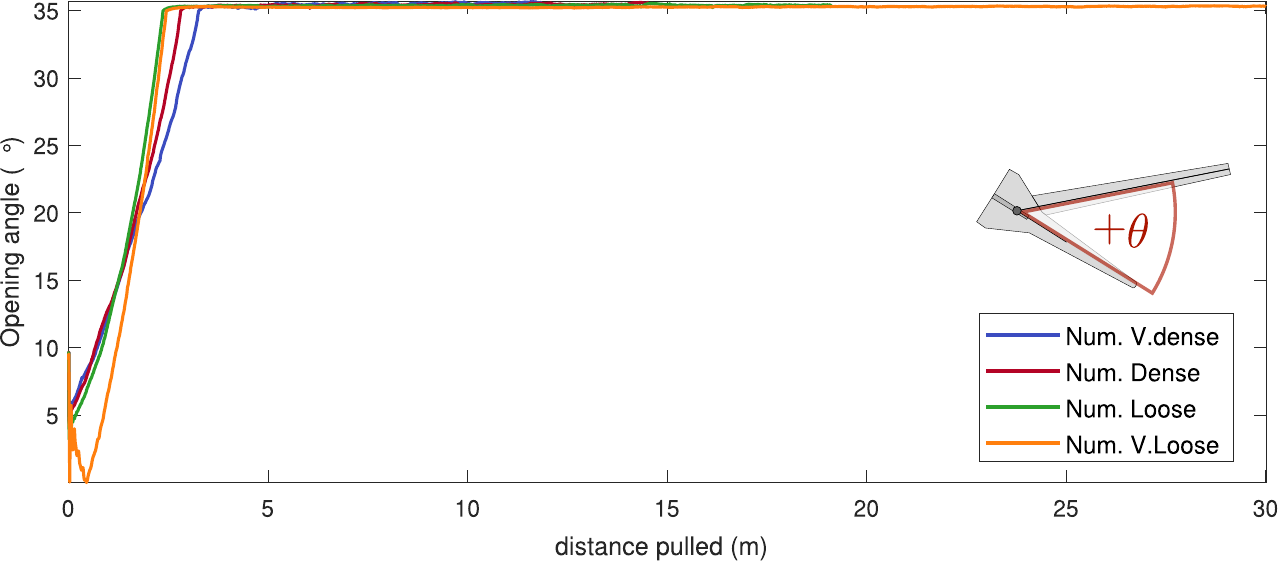}
        \caption{Opening angle.}
        \label{fig: anchor_opening_angle}
    \end{subfigure}
   
    \begin{subfigure}[b]{0.8\linewidth}
        \centering
        \includegraphics[width=\linewidth]{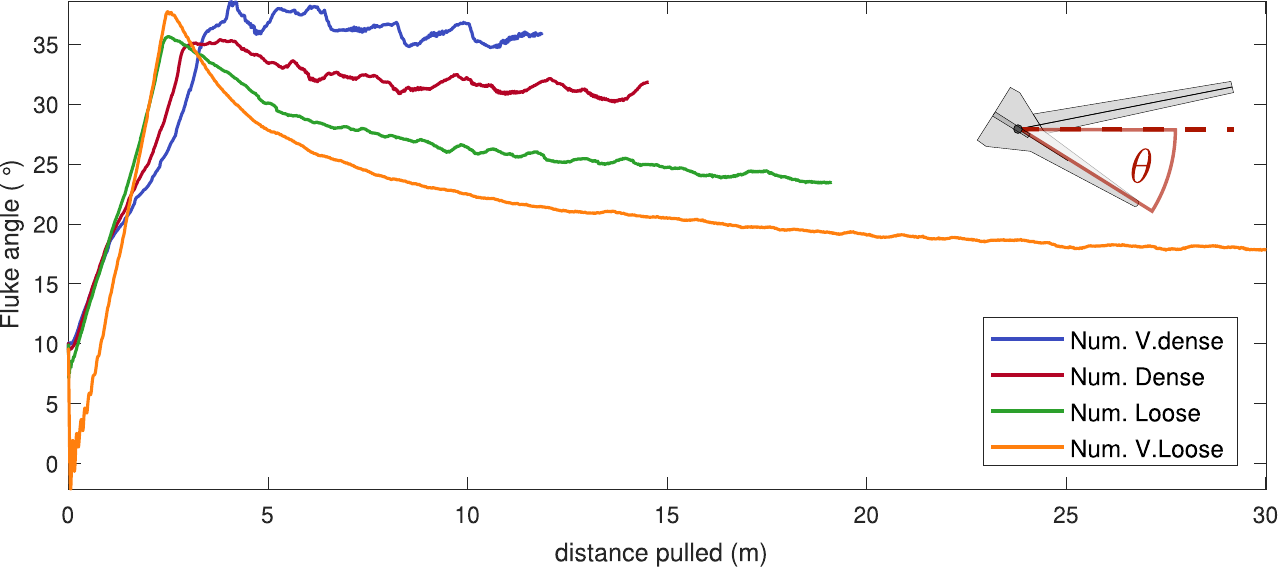}
        \caption{Fluke angle to the horizontal.}
        \label{fig: anchor_fluke_angle}
    \end{subfigure}
    
    \begin{subfigure}[b]{0.8\linewidth}
        \centering
        \includegraphics[width=\linewidth]{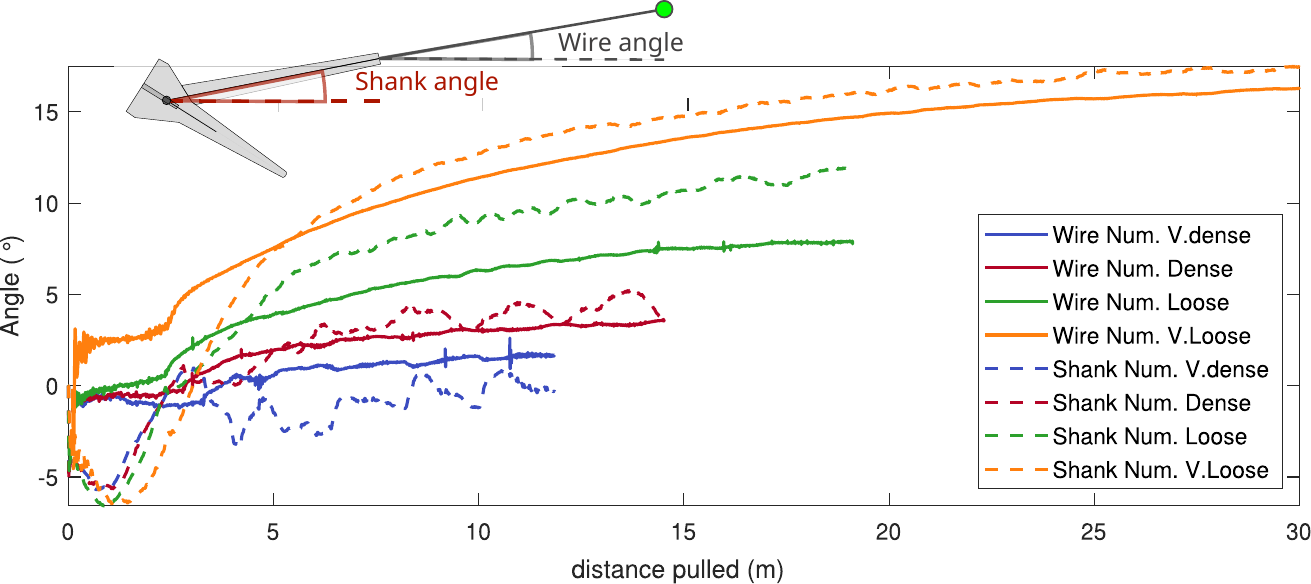}
        \caption{Wire and shank angle to the horizontal.}
        \label{fig: anchor_shank_angle}
    \end{subfigure}
    \caption{AC-14 anchor penetration: Opening angle (a) and, fluke (b) and, wire and shank angle (c), both to horizontal for the AC-14 anchor with drag distance.}
    \label{}
\end{figure}

The kinematic behaviour of the anchor is very different in very dense and loose sand conditions. This can be seen by observing the differences between Figure \ref{fig:rd_82_full_plot.pdf}  and Figure \ref{fig:rd_4_full_plot.pdf}, where the latter shows the anchor position in the very loose sand after $17$m of drag. In the very loose case, the shank and fluke are inclined at $15^\circ$ and $20^\circ$ to the horizontal, respectively. The result for opening angle and, fluke and shank angle, both to horizontal with drag distance are shown respectively in Figures \ref{fig: anchor_opening_angle}, \ref{fig: anchor_fluke_angle} and \ref{fig: anchor_shank_angle}. For the four densities the only common response is the opening angle in Figure \ref{fig: anchor_opening_angle} which all quickly reach the maximum, with the lowest density first, of $35^\circ$. For the fluke and shank angle, Figures \ref{fig: anchor_fluke_angle} and \ref{fig: anchor_shank_angle} show that for all densities a convergence in the angle is achieved with drag distance however the behaviour of the anchor is different. For the very dense sample the shank angle reduces towards zero with drag distance, whereas when the density is reduced the angle increases with distance, additionally the ultimate value shank angle also increases with reducing relative density. The opposite is true for the fluke angle, where for the lower densities the fluke is becoming more parallel to the horizontal. However, for the high densities the fluke angle is largest and changes marginally, with the majority of the change being oscillatory. 

In general the response of the anchor for the range of sands considered follows the description in Chapter 3 of \cite{puech1984use} where a number of full scale experiments were performed on articulated (shank and fluke are free to rotate) anchors. In summary the shank angle stays approximately parallel to the wire/chain, see Figure \ref{fig: anchor_shank_angle}, which is pulling the anchor whilst the fluke is initially angled down in the soil causing penetration. For the loose samples 
as the anchor is pulled the angle of the fluke creates a resultant forces which pushes the anchor deeper into the sample, the anchor subsequently rotates so that the fluke becomes more parallel to the horizontal whist the shank becomes more vertical. Whereas for the denser samples the majority of the penetration is caused by the fluke angle. This is also noted in \cite{puech1984use}, particularly when the soil is very dense and the anchor is "tripped" causing the anchor to fully open with only the fluke tip embedding.

\subsection{Influence of anchor size}\label{sec: anchor size}

\subsubsection*{Setup and scope}

The scope of this numerical study is to explore the impact of anchor size (or mass) on the penetration depth in uniform sands. Two sand densities are considered (loose and very dense), the length dimensions of the 8.7 tonne AC-14 anchor, used in the rest of the paper, are scaled by a factor of 0.5, 1.5 and 2.0 times their original values. Overall the anchors were dragged $6$ times their fluke length with their dimensions, and initial penetration due to self weight, shown in Table \ref{tab:anchorSizes}.

All simulations used the same mesh as in Section \ref{sec: anchor prediction} but with their element size and overall dimensions scaled by the corresponding factor in Table \ref{tab:anchorSizes}, the pull velocity is also scaled by the same factor. The only parameter which was not scaled was the time step increment.

All depth and pull distance results in this section are normalised with respect to the fluke length.

\subsubsection*{Results and discussion}
\begin{figure}[ht!]
    \centering
    \includegraphics[width=1.0\linewidth]{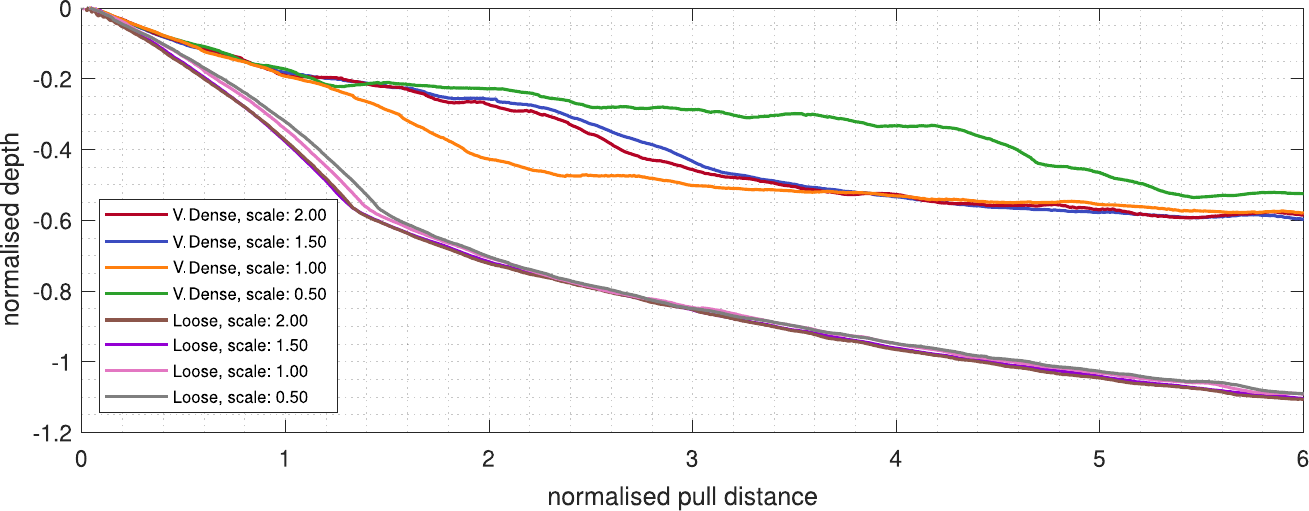}
    \caption{AC-14 scaled anchor penetration: scaled penetration depths for loose and very dense sands. Results are normalised with respect to the fluke length.}
    \label{fig:scale_results.pdf}
\end{figure}

The scaled anchor penetration results for the loose and very dense case are shown in Figure \ref{fig:scale_results.pdf}. The figure shows that regardless of scale, the normalised penetration depths achieve a very similar penetration, this is particularly true for the loose case where the penetration curves nearly all overlap exactly. For the very dense case the penetration of the fluke happens at different normalised drag distances. However, once the fluke has penetrated into the sand the normalised depths for all scales are very similar, see the normalised depth for the drag distance in the range $5.5$-$6.0$. The initial penetration of the anchors due to self weight is provided in Table \ref{tab:anchorSizes}, and shows that initial penetration does not scale proportional to the fluke length. An extreme example is the comparison between the $2$ and $0.5$ scaling factors in the loose sand condition, where the initial penetration due to self weight is $8.67$ times larger, despite the fluke only being $4$ times larger. Therefore when considering anchor penetration depths, the influence of the anchor size on the initial penetration should be taken into account as it has the potential to increase the overall penetration depth by 5-10\% depending on the sand state and the size of the anchor. It is interesting that for the scaling factors of $0.5$ and $1.0$, the normalised initial penetration is very similar and for both densities with $0.5$ having a larger normalised penetration, whilst the opposite trend is seen for the larger anchors. This is due to the non-planar geometry of the fluke which results in different anchors having different contact patch sizes with the sand which deforms to support the weight of the anchor.  

\begin{table}[ht!]
\centering
\caption{AC-14 anchor fluke lengths, $l_f$, mass and initial embedment under self weight, $d_{\text{in}}$, in loose and very dense conditions.}
\label{tab:anchorSizes}
{\small \begin{tabular}{c|c|c|c|c|c|c}
       scaling &  $l_f$ & mass  &  loose $d_{\text{in}}$ &  v. dense $d_{\text{in}}$ & loose  & v. dense \\ 
        factor & (m) & (tonnes) & (m) & (m) & ~~$d_{\text{in}}/l_f$~~ & $d_{\text{in}}/l_f$ \\ \hline
      $0.5$ & $0.85$  & $1.09$  & $0.052$ & $0.028$ & $0.061$ & $0.033$\\
      $1.0$ & $1.70$  & $8.70$  & $0.102$ & $0.050$ & $0.059$ & $0.029$\\
      $1.5$ & $2.55$  & $29.36$ & $0.203$ & $0.130$ & $0.079$ & $0.051$\\
      $2.0$ & $3.40$  & $69.60$ & $0.451$ & $0.234$ & $0.131$ & $0.068$ \\ \hline 
\end{tabular}}
\end{table}

\begin{figure}[ht!]
    \centering
    \begin{subfigure}[t]{0.324\textwidth}
        \centering
       \includegraphics[width=\textwidth]{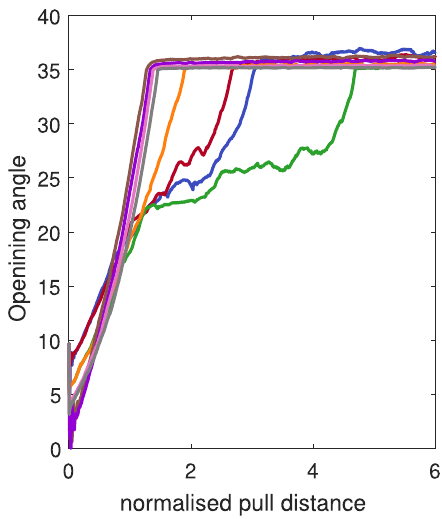}
        \caption{} \label{fig:scaled_opening_angle}
    \end{subfigure}%
    ~ 
    \begin{subfigure}[t]{0.324\textwidth}
        \centering
       \includegraphics[width=\textwidth]{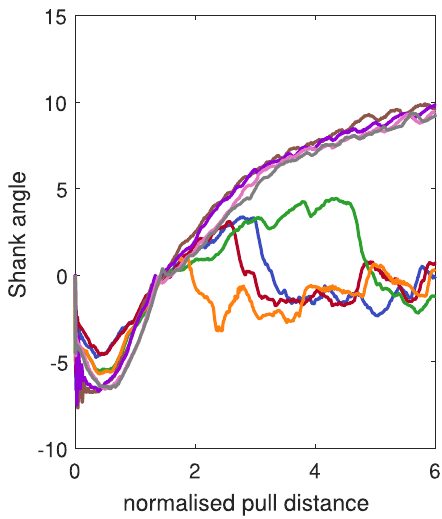}
        \caption{} \label{fig:scaled_shank_angle}
    \end{subfigure}%
    ~
     \begin{subfigure}[t]{0.324\textwidth}
        \centering
       \includegraphics[width=\textwidth]{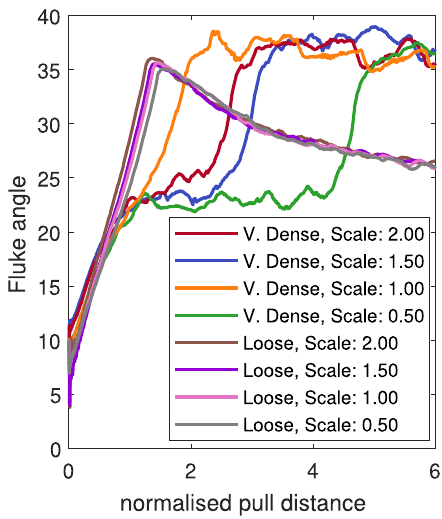}
        \caption{} \label{fig:scaled_fluke_angle}
    \end{subfigure}%
    \caption{AC-14 scaled anchor: (a) opening angle, (b) shank angle and (c) fluke angle against normalised pull distance. Results are normalised with respect to the fluke length.}\label{fig: scaled angles}
\end{figure}

A similar response is seen in the angles related to the anchor, see Figure \ref{fig: scaled angles}. For the loose case, the opening, shank and fluke angles follow almost the same trend. For the very dense case, the trends are not the same however their values at approximately the normalised pull distance of $5$, and further, are very similar. The response of the normalised depth for the very dense case is reflected in the change in the angles, for example at the scale $0.5$, the depth starts to increase substantially at a drag distance of $4.5$m (see Figure \ref{fig:scale_results.pdf}), caused by the increase in the opening and fluke angles at a drag distance of $4.5$m seen in Figure \ref{fig: scaled angles}. The delay in the normalised penetration of this scale is caused by the delay in the rotation of the anchor, but when the anchor starts to rotate, and therefore penetrate the seabed, the overall response is very similar. 

\begin{figure}[ht!]
        \centering
       \includegraphics[width=0.6\textwidth]{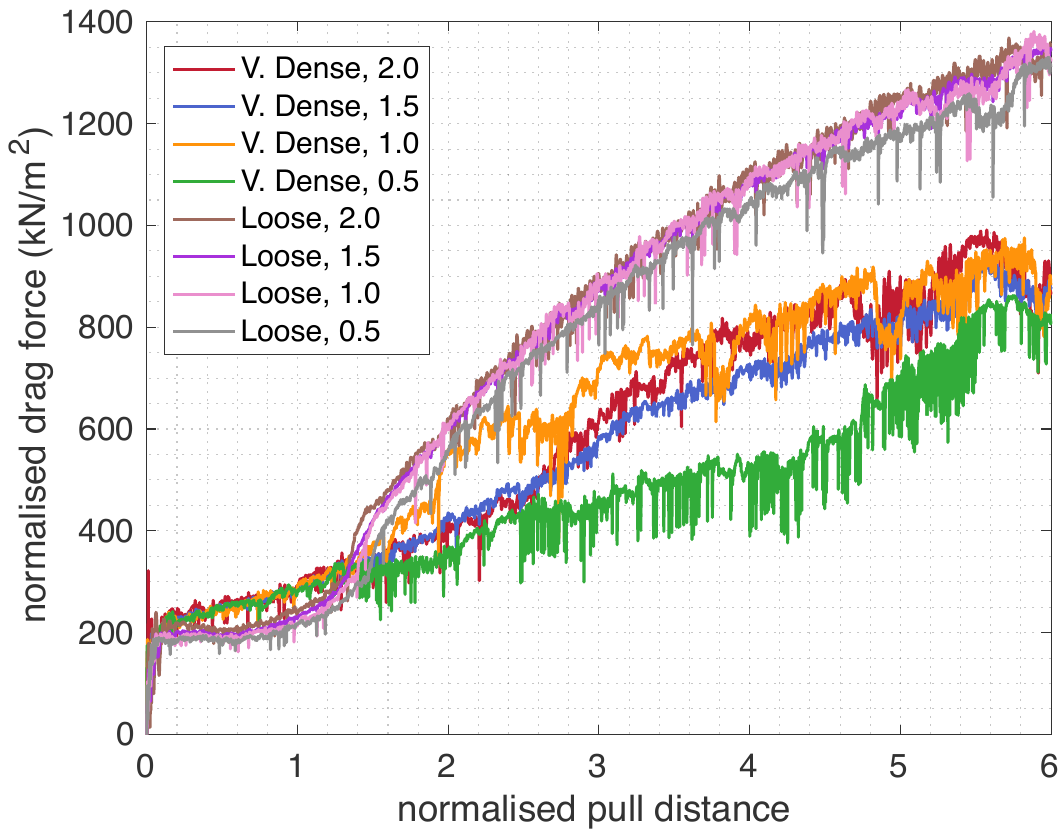}
    \caption{AC-14 scaled anchor: normalised drag force versus normalised drag distance. }
    \label{fig:scaled_force}
\end{figure}

Figure~\ref{fig:scaled_force} provides the drag force normalised by the length of the fluke squared (i.e. the relative change of the anchor area) versus the drag distance normalised by the fluke length for the four anchor sizes in loose and very dense sand. The consistency of the loose results demonstrates that the drag force scales with the area of the anchor and although there are variations in the very dense results, which are consistent with penetration depth and angle results discussed previously, the final normalised force is common across the anchor sizes.    

\subsection{Influence of anchor geometry}

The previous analyses in this paper focused on an AC-14 anchor. This section investigates the impact of changing the anchor geometry by comparing the AC-14 results with the penetration depths obtained from a 7.5 tonne Hall's stockless anchor (named after John Francis Hall, 1854-1897), another widely used shipping anchor. All the mesh and material parameter values are the same as in Section \ref{sec: anchor prediction}, the coefficient of friction acting at the anchor-sand interface is $0.3$. The penetration results of the Hall anchor are compared to experiments obtained using the same experimental setup as Sharif ~\emph{et al.} \cite{Yaseen2024}.

\subsubsection*{Setup and scope}
\begin{figure}[ht!]
    \centering
    \begin{subfigure}[t]{0.49\textwidth}
        \centering
       \includegraphics[width=1\textwidth]{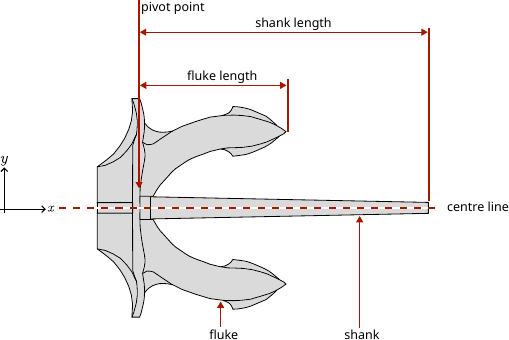}
        \caption{} \label{fig: hall anchor_schematic1}
    \end{subfigure}%
    \begin{subfigure}[t]{0.49\textwidth}
        \centering
      \includegraphics[width=1\textwidth]{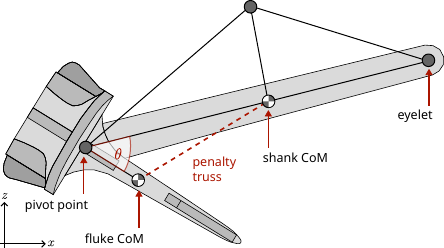}
        \caption{}\label{fig: hall anchor_schematic2}
    \end{subfigure}

    \begin{subfigure}[t]{1\textwidth}
        \centering
      \includegraphics[width=1\textwidth]{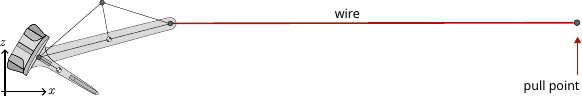}
        \caption{}\label{fig: hall anchor_schematic3}
    \end{subfigure}
    \caption{Hall anchor: (a) top down schematic, (b)side profile and (c) extended profile. The maximum opening angle is $45^\circ$.}\label{fig: hall anchor_schematic}
\end{figure}

The Hall anchor is dragged though 4 homogeneous sands, defined in the CPT calibration, Section \ref{sec: CPT calibration}, as \textit{very dense},  \textit{dense}, \textit{loose} and \textit{very loose}. The schematic of the anchor pull is shown in Figures \ref{fig: hall anchor_schematic1}, \ref{fig: hall anchor_schematic2} and \ref{fig: hall anchor_schematic3}, where the wire for this simulation has a length of $9$ m (different to the AC-14 at $9.5$ m). Note that the maximum opening angle of the Hall anchor is $45^\circ$.

The corresponding total mass, rotational inertia (rot. inertia) and geometric parameters (corresponding to Figure \ref{fig: hall anchor_schematic}) are shown in Table \ref{tab: Hall anchor mass, intertia and length properties}.
\begin{table}[ht!]
\centering
\caption{Hall anchor: Total mass, rotational inertia (rot. inertia) and geometric properties of the anchor components, where the CoM position and length are measured from the pivot point.}
\label{tab: Hall anchor mass, intertia and length properties}
{\small \begin{tabular}{r|c|c|c|c}
      & mass (kg)  & rot. inertia (kgm$^2$) & CoM position (m) & length (m) \\ \hline
fluke & 5688.4 & 1376               & 0.232                   &  1.75 \\ 
shank & 1811.6 & 1643               & 1.328                   & 3.23 \\ \hline 
\end{tabular}}
\end{table}
Finally the pull point position is set to coincide with the initial surface of the soil ($z=8$m in all of the examples in this paper) and has a prescribed velocity of $\{v\}=\{0.1 \quad 0 \quad 0\}^T~\text{m/s}$.
\subsubsection*{Results and discussion}

\begin{figure}
    \centering
    \includegraphics[width=\textwidth]{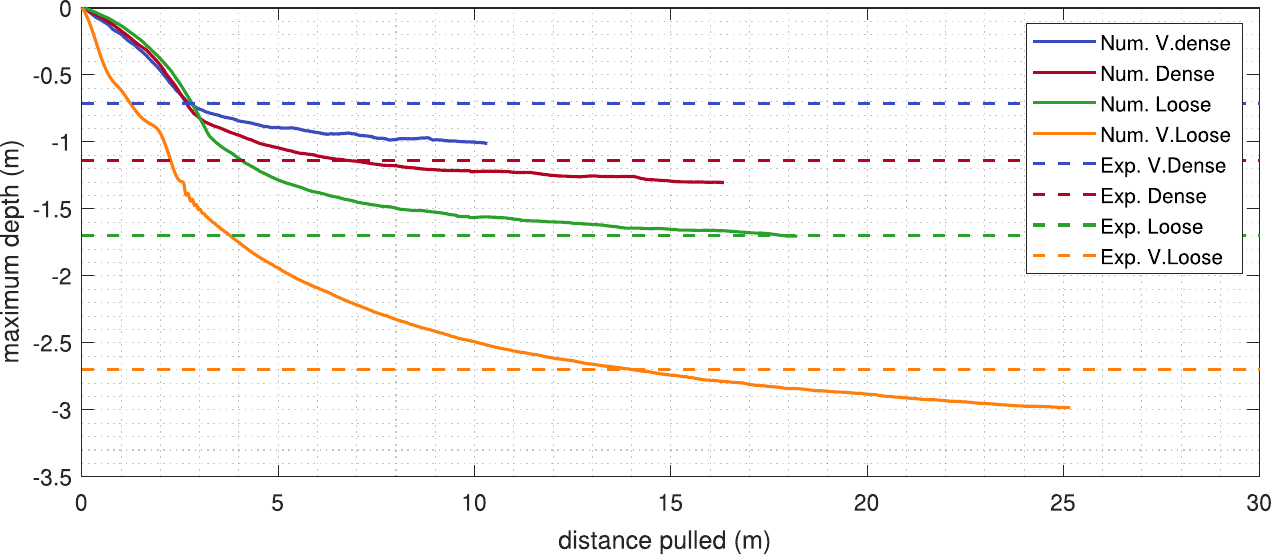}
    \caption{Hall anchor penetration: Comparison of anchor trajectories with the final depth of the experimental results.}
    \label{fig:Hall Anchor penetration}
\end{figure}

\begin{figure}[ht!]
    \centering
    \begin{subfigure}[t]{0.324\textwidth}
        \centering
       \includegraphics[width=\textwidth]{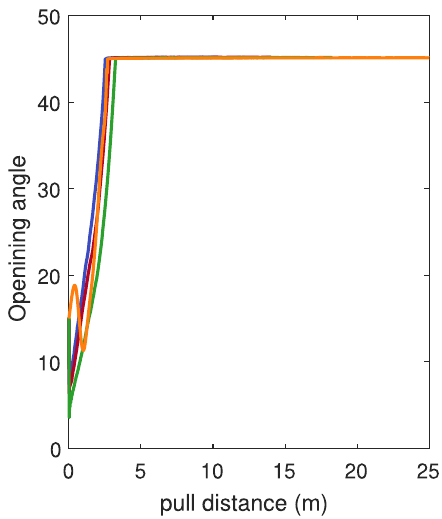}
        \caption{} \label{fig:hall_opening_angle}
    \end{subfigure}%
    ~ 
    \begin{subfigure}[t]{0.324\textwidth}
        \centering
       \includegraphics[width=\textwidth]{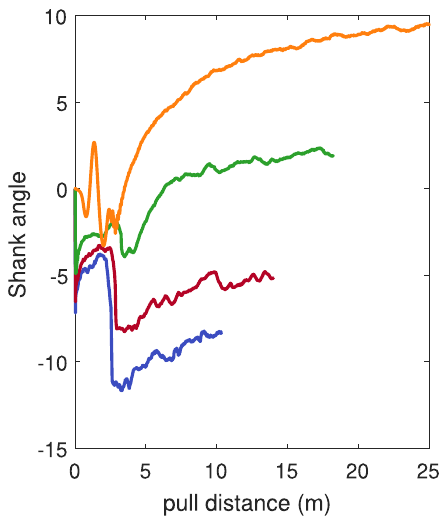}
        \caption{} \label{fig:hall_shank_angle}
    \end{subfigure}%
    ~
     \begin{subfigure}[t]{0.324\textwidth}
        \centering
       \includegraphics[width=\textwidth]{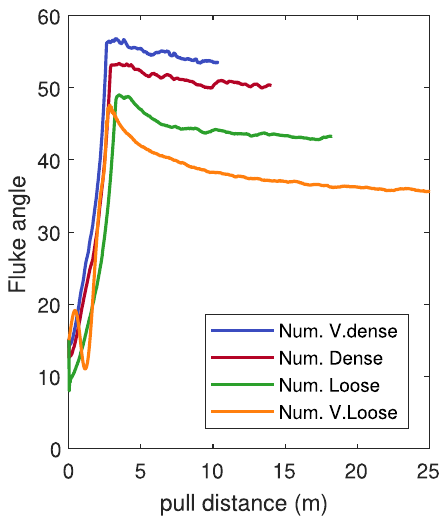}
        \caption{} \label{fig:hall_fluke_angle}
    \end{subfigure}%
    \caption{Hall anchor: The opening (a), shank (b) and fluke angle (c) against pull distance.}\label{fig: hall angles}
\end{figure}

The numerically predicted Hall anchor penetration depths are show in Figure \ref{fig:Hall Anchor penetration}, with the equivalent experimental data also shown for the final penetration depth. The opening, shank-horizontal and fluke-horizontal angles are also shown respectively in Figures \ref{fig:hall_opening_angle}, \ref{fig:hall_shank_angle} and \ref{fig:hall_fluke_angle}.

\begin{figure}[ht!]
    \centering
    \begin{subfigure}{0.55\linewidth}
        \centering
        \includegraphics[width=\linewidth]{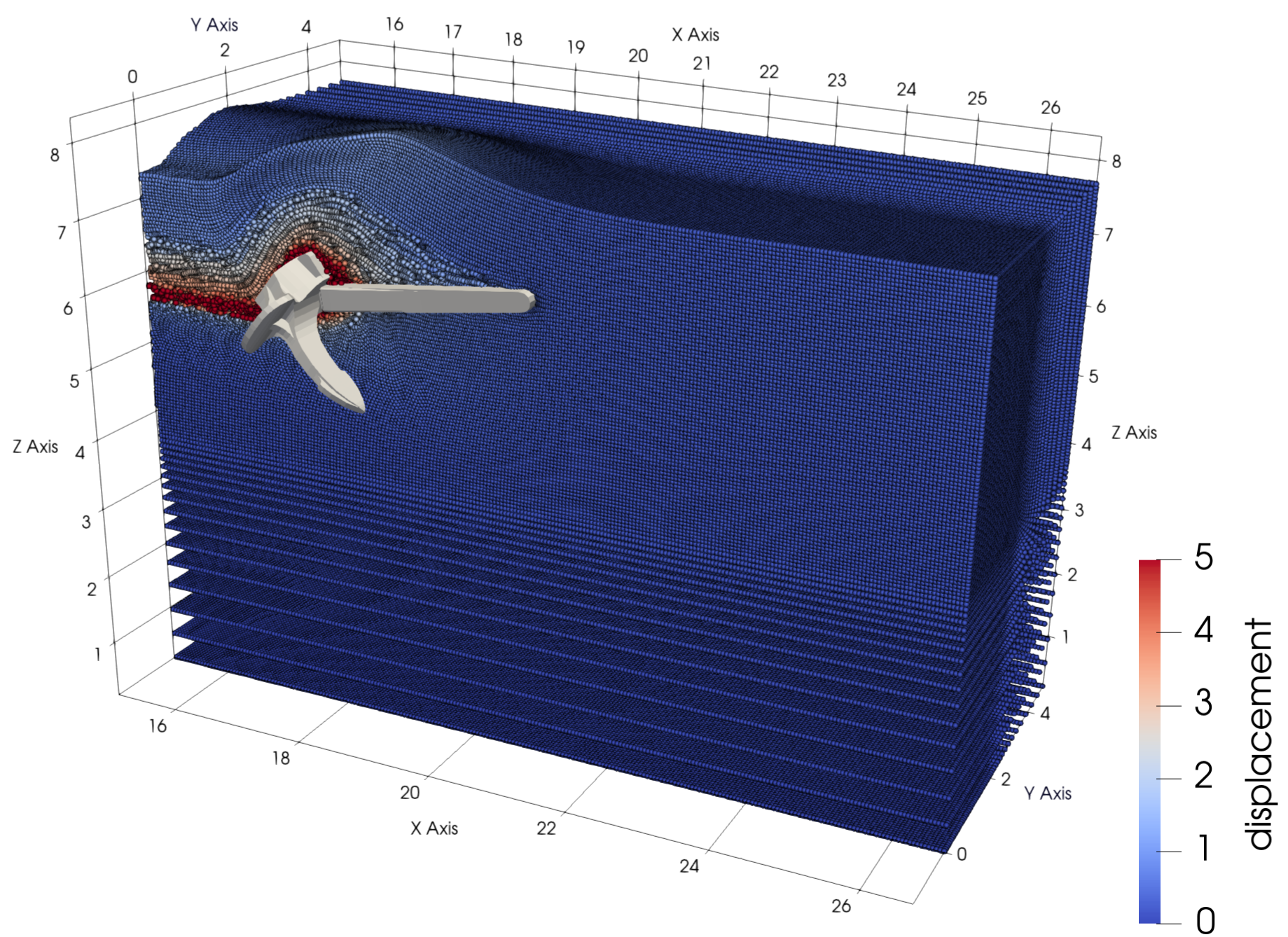}
        \caption{}
        \label{fig:hall_loose}
    \end{subfigure}%
    \hfill
    \begin{subfigure}{0.44\linewidth}
        \centering
        \includegraphics[width=\linewidth]{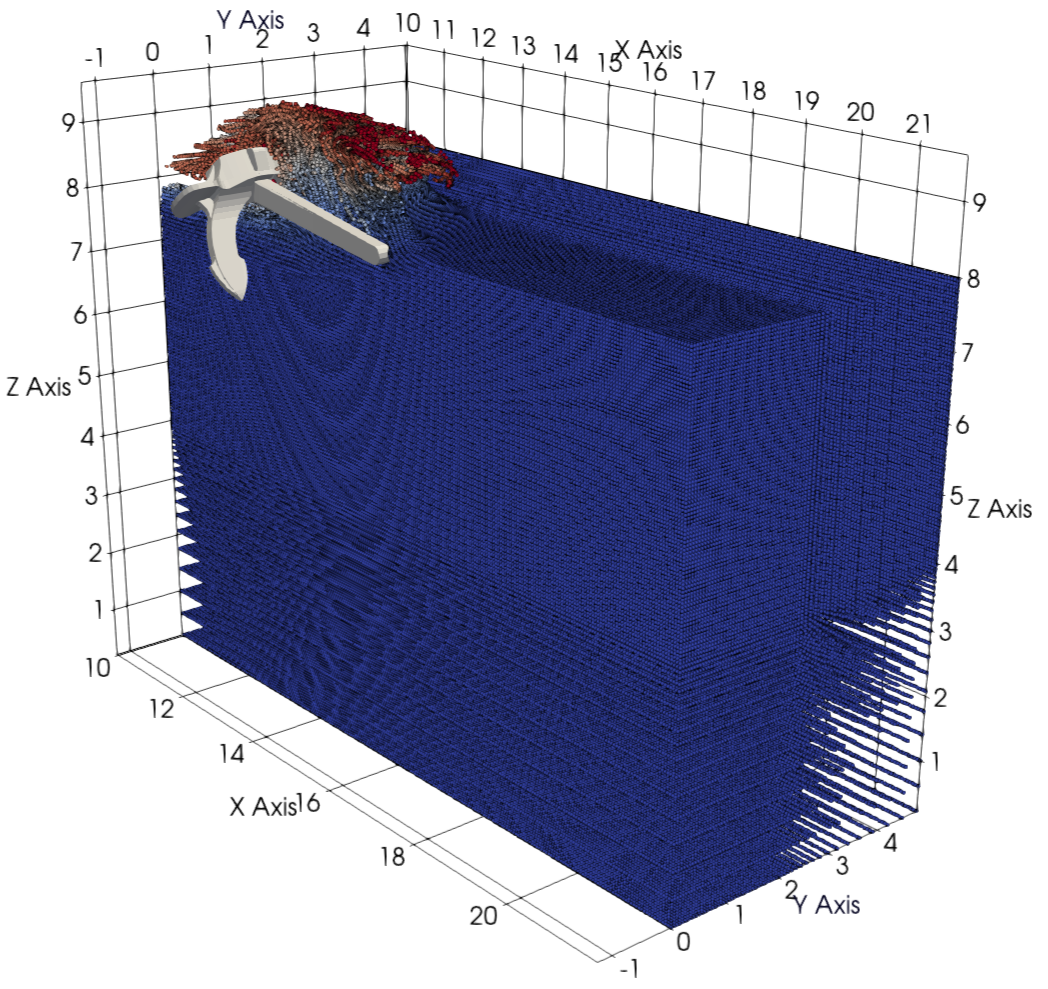}
        \caption{}
        \label{fig:hall_dense}
    \end{subfigure}
    \caption{Hall anchor: anchor position for the very loose at a drag distance of $15$m (a) and very dense sands at a drag distance of $9$m (b).}
    \label{fig:hall_anchor_positions}
\end{figure}

The comparison of the numerical and experimental results of the Hall anchor are consistent with the AC-14 anchor. For the very dense numerical result the depth is over predicted, this is due to nearly all the penetration from the fluke rotating into the sand, Figure~\ref{fig:hall_fluke_angle}. Very good agreement is achieved for the other sand densities, with the very loose result slightly over predicting the penetration at $3$m. For a qualitative comparison see Figure~\ref{fig:hall_anchor_positions} which shows the very loose and dense anchor position in the sand at their steady state.
The validation of the penetration depth for the Hall anchor further support the predictive capability of the proposed methodology, it shows the method is robust to a large range of sand relative densities and different anchor designs.

\section{Discussion and practical implications}

\subsection{Anchor performance}

\begin{figure}[ht!]
        \centering
       \includegraphics[width=0.6\textwidth]{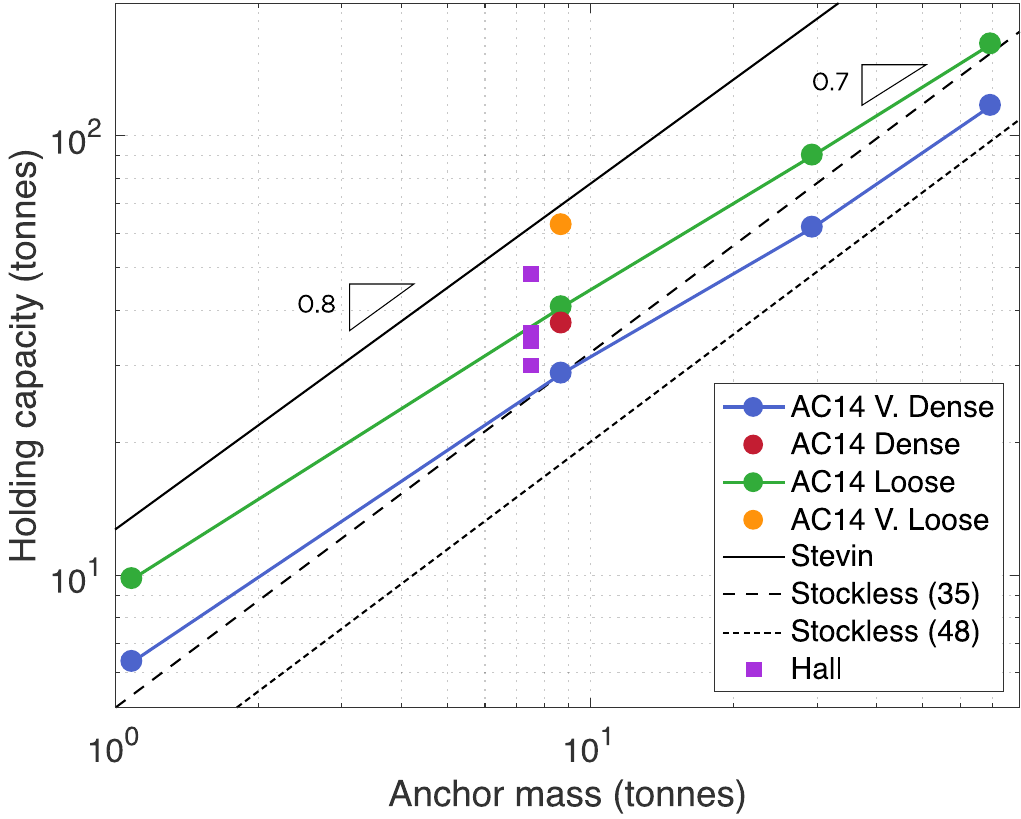}
    \caption{Anchor performance: holding capacity with anchor mass. AC-14 and Hall numerical results compared with Stevin\textsuperscript{\textregistered} and US Navy Stockless ($35^\circ$ and $48^\circ$ fluke angles) lines from \cite{NCL1987}. }
    \label{fig:scaled_capacity}
\end{figure}

It was shown in Section~\ref{sec: anchor size} that the capacity of a drag anchor scales with the square of the relative change in the anchor fluke length (i.e. linearly with the relative change of the anchor area). However, it is more typical to express the holding capacity of an anchor in terms of anchor mass (see for example \cite{vryhof2018manual}). Figure~\ref{fig:scaled_capacity} provides the holding capacity for different masses of AC-14 anchor in loose and very dense sand conditions. In both cases the holding capacity scales with the mass of the anchor to the power of $0.7$. This value is slightly lower than that suggested for Vryhof Stevin\textsuperscript{\textregistered} and US Navy Stockless anchors in \cite{NCL1987} who suggest a power of $0.8$ based on field trials. However, these field trials bound the numerical modelling result and provide additional validation of the MPM simulations presented in this paper. The numerical results also go some way in addressing the point raised by Randolph and Gourvenec \cite{Randolph2011} that design curves, such as those presented in Figure~\ref{fig:scaled_capacity}, have been historically based on \emph{``simple extrapolation of the ultimate holding capacity measured in small-scale tests to larger anchors.’’} The results confirm that the sand relative density shifts the position of the ultimate capacity, which is intrinsically linked to the anchor penetration depth, which is dependent on the sand state. This is an important point as it is typical that design charts combine sand and hard clay into a single performance line, often with separate lines for medium and soft clay. The results in Figure~\ref{fig:scaled_capacity} suggest that this may not be appropriate for all anchor designs.  However, changing the relative density does not appear to change the exponent linking the capacity to the anchor mass. Discrete points are also included on Figure~\ref{fig:scaled_capacity} for the 7.5 tonne Hall anchor in very loose to very dense sand conditions and the 8.7 tonne AC-14 anchor in very loose and dense sand conditions. These discrete points confirm the trends of very dense and loose AC-14 data and fall within the field trial data \cite{NCL1987}.

\begin{figure}
    \begin{subfigure}[t]{0.49\textwidth}
    \centering
    \includegraphics[width=\linewidth]{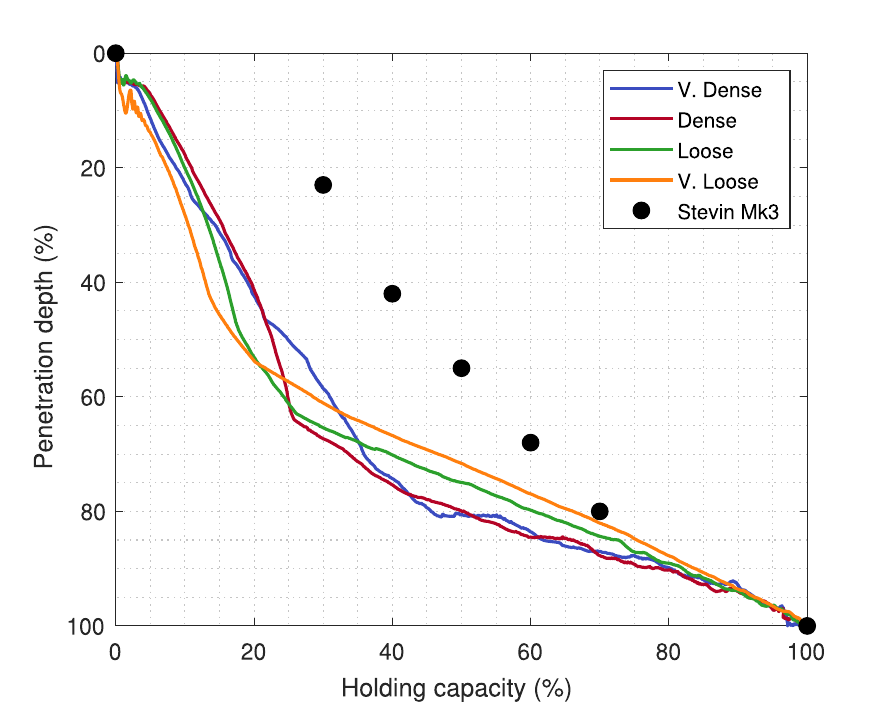}
    \caption{AC-14}
    \label{fig:AC14_HC}
    \end{subfigure}\hfill 
    \begin{subfigure}[t]{0.49\textwidth}
    \centering
    \includegraphics[width=\linewidth]{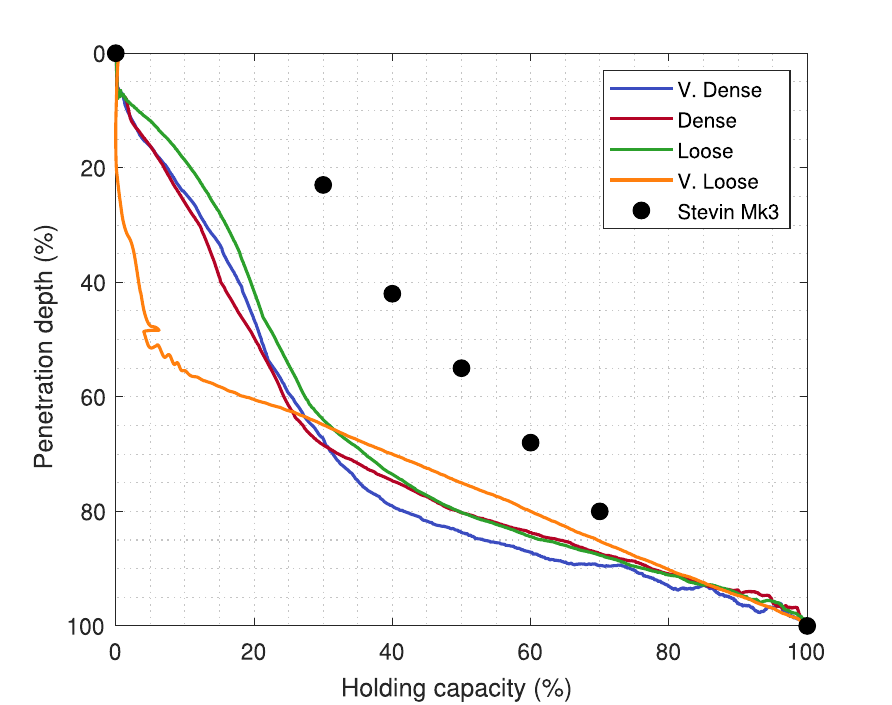}
    \caption{Hall}
    \label{fig:Hall_HC}
    \end{subfigure}
    \caption{Anchor performance: percentage penetration depth with holding capacity in different sand conditions. Final absolute values for the penetration depth and holding capacity are provided in Table~\ref{tab:holding capacity}.}
    \label{fig:holding_capacity}
\end{figure}

Drag anchor performance is also often presented in the form shown in Figures~\ref{fig:AC14_HC} and \ref{fig:Hall_HC} showing  percentage of ultimate holding capacity versus percentage of ultimate penetration depth (in the figure for the AC-14 and Hall anchors in particular). The final absolute values for the penetration depths and holding capacities are provided in Table~\ref{tab:holding capacity}. The MPM results are included for very loose to very dense sand conditions and compared against data from Vryhof \cite{vryhof2018manual} for the Stevin\textsuperscript{\textregistered} Mk3. The response of the AC-14 and Hall anchors are similar except in the case of very loose sand where the Hall anchor requires significantly more relative penetration to generate an equivalent percentage of its holding capacity. However, the shape of the AC-14 and Hall anchor curves are different to the Stevin\textsuperscript{\textregistered} Mk3 anchor, which generates more of its percentage capacity with a lower proportion of its penetration depth. This highlights the importance of anchor-specific performance charts versus utilisation of \emph{standard} curves.

\begin{table}[ht!]
\centering
\caption{Ultimate Holding Capacity (UHC) and maximum penetration, $d_p$, for the AC-14 and Hall anchors. Initial penetration under self weight is not included.}\label{tab:holding capacity}

{\small \begin{tabular}{lcc | cc}
& \multicolumn{2}{c|}{AC-14} & \multicolumn{2}{c}{Hall}\\
sand state & UHC (kN) &  $d_p$ (m) & UHC (kN) &  $d_p$ (m)\\
\hline
v. dense  & 283.5 & 1.03 & 294.3 & 1.01\\ 
dense     & 360.0 & 1.38 & 341.0 & 1.30\\ 
loose     & 491.5 & 2.21 & 323.7 & 1.70\\ %
v. loose  & 601.4 & 3.51 & 485.7 & 2.98 \\  
\hline
\end{tabular}}
\end{table}

\subsection{CBRA implications}

The UK Carbon Trust's Cable Burial Risk Assessment (CBRA, \cite{CBRA}) assumes that the penetration depth of a drag anchor can be approximated as 
    \begin{equation}\label{eqn:CBRA}
        d_p= S_f \sin(\theta_{\text{sf}})  l_f,
    \end{equation}
where $S_f$ is the seabed factor (taken as $1$ for sands in the CBRA guidance), $\theta_{\text{sf}}$ is the shank-fluke maximum opening angle (or $45^\circ$ if the maximum opening angle is not known) and $l_f$ is the fluke length. The results presented in this paper highlight two key shortcomings of the UK Carbon Trust's CBRA assumptions in terms of anchor penetration:
\begin{enumerate}
    \item  Equation (\ref{eqn:CBRA}) assumes that the shank is horizontal at the seabed surface and the anchor is fixed at its maximum opening angle. The numerical results have shown that the shank is only horizontal for the very dense sands. In other cases the shank-fluke joint is below the seabed surface and the shank inclined to the horizontal. This challenges the geometric assumptions in Equation (\ref{eqn:CBRA}).
    \item The CBRA guidance document gives a constant seabed factor of $S_f=1$ for sands, suggesting that an anchor will achieve the same penetration depth irrespective of the sand relative density. This is incorrect based on the numerical results of this section and the physical testing of Sharif \emph{et al.} \cite{Yaseen2024}. This constant seabed factor is linked to Point 1, where it is assumed that the shank is horizontal and at the seabed surface through the anchor drag process. 
\end{enumerate}
The anchor scaling tests have shown that the penetration depths associated with the anchor drag process vary linearly with fluke length, $l_f$, confirming the dependency in Equation (\ref{eqn:CBRA}). These tests have also highlighted that initial penetration of the anchor under self weight is not linearly dependent on the fluke length; this initial self weight penetration is 5-10\% of the penetration achieved during the drag process over the tested range of scaling factors.  

Whilst acknowledging the shortcoming of the CBRA equation, the MPM anchor penetration results represent a chance to provide new sand density dependent seabed penetration factors.   Table~\ref{tab:CBRA} provides equivalent CBRA seabed factors for the AC-14 anchor calculated from
\begin{equation}
    \bar{S}_f^{(\theta_{sf})} = \frac{d_p}{l_p~\sin(\theta_{sf})},
\end{equation}
where $d_p$ is taken as the penetration depth obtained from the MPM analyses during the drag process, based on a fluke length of $l_f=1.7$m (i.e. the results shown in Figure~\ref{fig:Anchor penetration}). Values are tabulated using the actual maximum opening angle for the anchor in the MPM simulations ($35^\circ$ in all cases) and $45^\circ$, as the default value in the CBRA guidance document. The back-calculated seabed factors for an AC-14 anchor vary between $1.1$ and $3.6$ for the $35^\circ$ opening angle and $0.9$ and $2.9$ for the $45^\circ$ opening angle, which are in good agreement with the physical testing \cite{Yaseen2024} derived seabed factors, which have a range of $0.6$ to $3.4$ assuming a $45^\circ$ opening angle. 

Table~\ref{tab:CBRA_Hall} provides equivalent CBRA seabed factors for the Hall anchor. In this case the opening angle of the anchor from the MPM simulations was $45^\circ$ for all sand relative densities, the same as the CBRA default value. The back-calculated seabed factors for the Hall anchor vary between $0.8$ and $2.4$, which are significantly less than the AC-14 $35^\circ$ opening angle results. The penetration factors for the two anchors with a $45^\circ$ opening angle are similar in dense and very dense sands, but the Hall anchor has less penetration potential in loose sands,  likely due to the large plate behind the fluke-shank interface, which is not present on the AC-14 anchor (see Figures~\ref{Fig:AC14 anchor schematic} and \ref{fig: hall anchor_schematic}). The back-calculated seabed factors for the AC-14 and Hall anchors suggest that the seabed factor of $S_f=1$ from the CBRA guidance document is non-conservative beyond the very dense case. Note that the initial penetration under self weight is not included; in this case the seabed factors would increase by $\approx 5$\% if taken into account. It should also be noted that these seabed factors are considered to be the worst case scenario as they are based on fully drained conditions. The work of Brown \emph{et al.} \cite{Brown2026rate} has shown that rate effects can reduce the penetration depth, especially in the case of loose sands. However, the exact drag speed of an anchor when it crosses a subsea cable would not be known \emph{a priori} when undertaking an CBRA and providing worst case seabed factors facilitates a conservative approach.   

\begin{table}[ht!]
\centering
\caption{CBRA AC-14 seabed factors for different sands based on anchor penetration, $d_p$, from the MPM analyses with an 8.7 tonne anchor with a fluke length of $l_f=1.7$m. Initial penetration under self weight is not included. }
\label{tab:CBRA}
\begin{tabular}{l|c|c|c|c}
       sand state & ~~$d_p$ (m)~~ & ~~$d_p/l_f$~~ & $\bar{S}^{(35^\circ)}_f$ & $\bar{S}^{(45^\circ)}_f$ \\ \hline
      v. dense & $1.03$ & $0.61$ & $1.06$ & $0.86$ \\
      dense    & $1.38$ & $0.81$ & $1.42$ & $1.15$ \\
      loose    & $2.21$ & $1.30$ & $2.27$ & $1.84$ \\
      v. loose & $3.51$ & $2.06$ & $3.60$ & $2.92$ \\ \hline
\end{tabular}
\end{table}


\begin{table}[ht!]
\centering
\caption{CBRA Hall seabed factors for different sands based on anchor penetration, $d_p$, from the MPM analyses with an $7.5$ tonne anchor with a fluke length of $l_f=1.75$m. Initial penetration under self weight is not included.}
\label{tab:CBRA_Hall}
\begin{tabular}{l|c|c|c}
       sand state & ~~$d_p$ (m)~~ & ~~$d_p/l_f$~~ & $\bar{S}^{(45^\circ)}_f$ \\ \hline
      v. dense & $1.01$ & $0.57$ & $0.82$  \\
      dense    & $1.30$ & $0.74$ & $1.05$  \\
      loose    & $1.70$ & $0.97$ & $1.37$  \\
      v. loose & $2.98$ & $1.70$ & $2.41$  \\ \hline
\end{tabular}
\end{table}

\section{Conclusions}
The paper extends the work of Bird \emph{et al.} \cite{bird2024dynamic} to enable computationally tractable simulations of geotechnical problems where a rigid body, formed of multiple parts, is pulled/pushed a long distance through a highly deformable domain. A partitioned domain approach is developed which makes the computational cost per load step independent of the total pull/push distance. The method was applied to the GIMPM by moving the domain boundaries closer to the rigid body and only defining the GIMPs within the boundaries as active. The method was validated for an anchor pull problem, it was shown that for a large range in material properties if the partitioned domain was sufficiently large the anchor trajectory is uneffected by the reduced domain size. From the validation it was found that it is only necessary for the domain length in front of the anchor to be twice the bounding box length in the pull direction, where the box contains the entire anchor. This created a considerable reduction in problem size as in some cases it was necessary to drag the anchor at least 15 times its drag length. 

An accurate predictive CPT-based anchor penetration framework was also presented. It has two steps: (i) determine the material properties of the sand bed by comparing numerical CPTs with a experimental/field data; (ii) use the partitioned domain approach to calculate the anchor trajectory and the maximum penetration depth. The framework was validated for a large range of relative density sands against experimental results for two articulated anchor designs, it achieved accurate results for the final penetration depths versus both centrifuge experiments and field trial data. The method is predictive - no numerical, or material, parameter tuning is performed with the range of results achieved through the variation of only four mechanical parameters in the soil constitutive models (two elastic, two plastic) which are determined from the corresponding CPT analysis. Variations in the anchor size in the numerical analysis has confirmed that the penetration depth and the holding capacity scale with the fluke length and the fluke length squared, respectively, and that the holding capacity depends on the sand state as well as the anchor mass. The simulations have also revealed two key issues with the CBRA anchor penetration estimation equation, namely the assumption: (i) of a constant seabed factor of $1$ irrespective of the sand state; and (ii) that the shank is horizontal and at the seabed surface.  Back analysis of the simulation results suggest seabed penetration factors between $0.8$ and $3.6$ depending on the sand state, assumed opening angle and the anchor type, highlighting the potentially non-conservatism of the CBRA suggested penetration factor for sands.   In summary, this paper provides a predictive framework for predicting anchor penetration, which has been shown to be accurate over a very wide range of sand relative densities. The underlying material model can be calibrated using site investigation data that is routinely available along a cabling route, which opens to door to perform site specific analysis for anchor penetration within the CBRAs framework.

\section*{Acknowledgements}
This work was supported by the Engineering and Physical Sciences Research Council [grant numbers EP/W000970/1, EP/W000997/1 and EP/W000954/1]. The second author was supported by funding from the Faculty of Science, Durham University.

\bibliographystyle{plain}
\bibliography{cas-refs.bib}

\appendix
\section{Determining \texorpdfstring{$\delta\theta$, $\Delta\ddot{\theta}$ and $\Delta\delta\theta$}{delta-theta, Delta ddot-theta and Delta delta-theta}}\label{App:Linearisation of theta}
In this section the definitions of $\delta\theta$, $\Delta\ddot{\theta}$ and $\Delta\delta\theta$ are presented. To make the derivation clearer the following abuse of notation, \textit{applied only in this section}, is used for the difference of position, velocity and acceleration for the truss nodes $A$ and $B$, see Figure \ref{fig:truss frame with rot}:
\begin{align}
    \bm{x}:=\bm{x}^{AB} & = \bm{x}_B - \bm{x}_A\\
    \bm{v}:=\bm{v}^{AB} & = \bm{v}_B - \bm{v}_A\\
    \dot{\bm{v}}:=\dot{\bm{v}}^{AB} & = \dot{\bm{v}}_B - \dot{\bm{v}}_A
\end{align}

A requirement is the definition of 
For the definitions of $\delta\theta$, $\Delta\ddot{\theta}$ and $\Delta\delta\theta$
For the variations and linearisation it is 
Before defining the linearisation of $\ddot{\theta}$, it is necessary to define the first and second derivatives of $\theta$, respectively
\begin{equation}\label{equ: div theta x}
    \frac{\partial\theta}{\partial\bm{x}} =
    \left[\begin{array}{c}
{-x_3}{/({x_1}^2 + {x_3}^2)} \\
{x_1}/{({x_1}^2 + {x_3}^2)}
\end{array}\right]
\end{equation}
and
\begin{equation}\label{equ: second div theta x}
    \frac{\partial^2\theta}{{\partial\bm{x}}^2} =
    \left[\begin{array}{cc}
A & B\\
B & A
\end{array}\right]
\quad\text{where}\quad 
\begin{array}{l}
A =(2{x_1}{x_3})/({x_1}^2 + {x_3}^2)^2 \\
B= -({x_1}^2 - {x_3}^2)/({x_1}^2 + {x_3}^2)^2.
\end{array}
\end{equation}
Additional the rotational velocity and acceleration have the respective forms
\begin{equation}
    \dot{\theta} = \frac{\partial\theta}{\partial t}= \frac{\partial\theta}{\partial \bm{x}}\cdot\frac{\partial \bm{x}}{\partial t}
\end{equation}
where $t$ is time and
\begin{align}
    \ddot{\theta} = \bm{v}\cdot\frac{\partial^2\theta}{\partial {\bm{x}}^2}\cdot\bm{v} +  \frac{\partial\theta}{\partial \bm{x}}\cdot\dot{\bm{v}}
\end{align}

\subsection{Variation of ${\theta}$}
As $\theta$ is a function of only $\bm{x}$ its variations takes the form,
\begin{equation}\label{equ:variation of theta}
    \begin{split}
        \delta\theta &= \frac{\partial \theta}{\partial \bm{x}}\cdot\frac{\partial\bm{x}}{\partial t} \delta t = \frac{\partial \theta}{\partial \bm{x}}\cdot\delta\bm{x}
    \end{split}
\end{equation}
where $ {\partial \theta}/{\partial \bm{x}}$ is defined in Equation \eqref{equ: div theta x}.

\subsection{Linearisation of $\delta{\theta}$}
The linearisation of $\delta{\theta}$, defined in Equation \eqref{equ:variation of theta}, is 
\begin{equation}
\begin{split}
    \Delta\delta\theta & = \Delta\left(\frac{\partial \theta}{\partial \bm{x}}\cdot\delta\bm{x}\right)  = \Delta\left(\frac{\partial \theta}{\partial \bm{x}}\right)\cdot\delta\bm{x}
\end{split}
\end{equation}
where
\begin{equation}\label{equ: linearisation of div theta x}
\begin{split}
    \Delta\left(\frac{\partial \theta}{\partial \bm{x}}\right) &= \frac{\partial^2\theta}{\partial \bm{x}^2} \cdot\frac{\partial \bm{x}}{\partial t}\Delta t = \frac{\partial^2\theta}{\partial \bm{x}^2} \cdot\Delta \bm{x}
    \end{split}
\end{equation}
and ${\partial^2\theta}/{\partial \bm{x}^2}$ is defined in Equation \eqref{equ: second div theta x}.

\subsection{Linearisation of \texorpdfstring{$\ddot{\theta}$}{ddot-theta}}
The linearisation of $\ddot{\theta}$ is defined $\Delta\ddot{\theta}$ and it takes the form
\begin{align}
    \Delta\ddot{\theta} &= \Delta\left(\bm{v}\cdot\frac{\partial^2\theta}{\partial {\bm{x}}^2}\cdot\bm{v}\right) + \Delta\left(\frac{\partial\theta}{\partial \bm{x}}\cdot\dot{\bm{v}}\right)\\
    & = \Delta\bm{v}\cdot\frac{\partial^2\theta}{\partial {\bm{x}}^2}\cdot\bm{v}
      + \bm{v}\cdot\frac{\partial^2\theta}{\partial {\bm{x}}^2}\cdot\Delta\bm{v}
      + \bm{v}\cdot\Delta\left(\frac{\partial^2\theta}{\partial {\bm{x}}^2}\right)\cdot\bm{v}\\
    & + \Delta\left(\frac{\partial\theta}{\partial \bm{x}}\right)\cdot\dot{\bm{v}} + \frac{\partial\theta}{\partial \bm{x}}\cdot\Delta\dot{\bm{v}}.
\end{align}
where $\Delta\left({\partial\theta}/{\partial \bm{x}}\right)$ is defined in Equation \eqref{equ: linearisation of div theta x} and the remaining linearisation of $\bm{v}$, $\dot{\bm{v}}$, and $\Delta\left({\partial^2\theta}/{\partial \bm{x}^2}\right)$ are provided below. The expressions for $\dot{\bm{v}}$ and $\bm{v}$, and their subsequent linearisation, take the Newmark form since the Newmark method is used to ingrate time. The linearisation of acceleration has the form
\begin{align}
    \dot{\bm{v}} = \dot{\bm{v}}_{n+1} &= \frac{\bm{x}_{n+1}-\bm{x}_{n}}{\beta ~\Delta_n t^2}-\frac{\bm{v}_n}{\beta~\Delta_n t}-\dot{\bm{v}}_n\left(\frac{1}{2\beta}-1\right)\\
   \Delta \dot{\bm{v}} = \Delta\dot{\bm{v}}_{n+1} &= \frac{\Delta\bm{x}_{n+1}}{\beta ~\Delta_n t^2} = \frac{\Delta\bm{x}}{\beta ~\Delta_n t^2}
\end{align}
where $\Delta_n$ is the time integration increment. The linearisation of the velocity is
\begin{align}
    \bm{v} = \bm{v}_{n+1} &= \gamma\frac{\bm{x}_{n+1}-\bm{x}_n}{\beta~\Delta_n t}+\bm{v}_n\left(1-\frac{\gamma}{\beta}\right)+\Delta_n t \dot{\bm{v}}_n\left(1-\frac{\gamma}{2\beta}\right);\\
    \Delta \bm{v} = \Delta\bm{v}_{n+1} &= \gamma\frac{ \Delta \bm{x}_{n+1}}{\beta ~\Delta_n t} = \gamma\frac{ \Delta\bm{x}}{\beta ~\Delta_n t}.
\end{align}
Last the linearisation of the term ${\partial^2\theta}/{\partial {\bm{x}}^2}$ is expressed as
\begin{equation}
\begin{split}
    \Delta\left(\frac{\partial^2\theta}{\partial {\bm{x}}^2}\right) &= \frac{\partial}{\partial\bm{x}}\left(\frac{\partial^2\theta}{\partial {\bm{x}}^2}\right)\cdot\Delta\bm{x} = \left[\begin{array}{cc}
    A & B\\
    C & D
    \end{array}\right]
    \quad\text{where}\\
    A &= \left[
  -\frac{2x_3(3x_1^2 - x_3^2)}{(x_1^2 + x_3^2)^3},\ 
  \frac{2x_1(x_1^2 - 3x_3^2)}{(x_1^2 + x_3^2)^3}
\right]^T \cdot \Delta\bm{x} \\
B &= \left[
  \frac{2x_1(x_1^2 - 3x_3^2)}{(x_1^2 + x_3^2)^3},\ 
  \frac{2x_3(3x_1^2 - x_3^2)}{(x_1^2 + x_3^2)^3}
\right]^T \cdot \Delta\bm{x} \\
C &= \left[
  \frac{2x_1(x_1^2 - 3x_3^2)}{(x_1^2 + x_3^2)^3},\ 
  \frac{2x_3(3x_1^2 - x_3^2)}{(x_1^2 + x_3^2)^3}
\right]^T \cdot \Delta\bm{x} \\
D &= \left[
  \frac{2x_3(3x_1^2 - x_3^2)}{(x_1^2 + x_3^2)^3},\ 
  -\frac{2x_1(x_1^2 - 3x_3^2)}{(x_1^2 + x_3^2)^3}
\right]^T \cdot \Delta\bm{x}
\end{split}
\end{equation}

\end{document}